\documentclass[11pt,a4paper]{article}

\usepackage[T2A]{fontenc}
\usepackage[cp1251]{inputenc}

\usepackage{amssymb,amsmath,amsthm,amsfonts,epsf,fancyhdr,epic,eepic,xspace,pstricks,textcomp}

\usepackage{epsfig,graphicx,rotating,array}

\ifx\LibLexLoaded\undefined\else\expandafter \fi
\let\LibLexLoaded\relax

\def\makeatletter{\catcode`\@=11\relax}

\begingroup \makeatletter
\gdef\atc@end{\catcode`\@\errmessage{Unmatched \string\restoreatcode}}
\gdef\atc@stack{\atc@end}
\gdef\atc@next#1\atc@end{\def\atc@stack{#1\atc@end}}
\gdef\saveatcode{%
 \expandafter\expandafter\expandafter\atc@next
 \expandafter\the\expandafter\catcode\expandafter`\expandafter\@
 \expandafter\space\expandafter\atc@next\atc@stack
}
\gdef\restoreatcode{\catcode`\@=\atc@stack}
\endgroup

\saveatcode\makeatletter

\def\tempcount{\lib@alloc0\count\countdef\insc@unt}
\def\tempdimen{\lib@alloc1\dimen\dimendef\insc@unt}
\def\tempskip{\lib@alloc2\skip\skipdef\insc@unt}
\def\tempmuskip{\lib@alloc3\muskip\muskipdef\@cclvi}
\def\tempbox{\lib@alloc4\box\chardef\insc@unt}
\def\temptoks{\lib@alloc5\toks\toksdef\@cclvi}
\def\tempread{\lib@alloc6\read\chardef\sixt@@n}
\def\tempwrite{\lib@alloc7\write\chardef\sixt@@n}

\def\lib@alloc#1#2#3#4#5{%
  \advance\count1#1by\@ne
  \ifnum\count1#1<#4\else\lib@alloc@err#2\fi
  #3#5=\count1#1\relax
}
\def\lib@alloc@err#1{\errmessage{No room for a new temporary #1}}

\let\lib@tmp\undefined
\let\lib@next\undefined

\def\lib@empty{} 

\let\lib@gtmp\undefined 

\newif\ifyes

\def\testfile#1{%
  \begingroup
  \tempread\lib@tmp
  \openin\lib@tmp=#1\relax
  \ifeof\lib@tmp
    \closein\lib@tmp \endgroup \yesfalse
  \else
    \closein\lib@tmp \endgroup \yestrue
  \fi
}

\restoreatcode

\saveatcode\makeatletter

\newtoks\evp@shadow 
\evp@shadow{\evp@hidden}
\def\evp@hidden{\evp@hidden}

\def\@atnextpar{}
\def\@ateverypar{}

\newif\if@evp 

\def\evp@setup{%
 \if@evp\else
  \@evptrue
  \evp@shadow\everypar
  \everypar{%
   \let\lib@tmp\@atnextpar \global\let\@atnextpar\lib@empty \lib@tmp
   \@ateverypar
   \the\evp@shadow
  }%
  \let\everypar\evp@shadow
  \aftergroup\evp@aftergroup
 \fi
}

\def\evp@aftergroup{%
 \ifx\@atnextpar\lib@empty
  \begingroup
   \edef\lib@tmp{\the\evp@shadow}%
   \ifx\lib@tmp\evp@hidden \else
    \global\everypar\evp@shadow
    \global\evp@shadow{\evp@hidden}%
   \fi
  \endgroup
 \else
  \evp@setup
 \fi
}

\def\atnextpar#1{%
 \evp@setup
 \begingroup
  \expandafter\toks@\expandafter{\@atnextpar #1}%
  \xdef\@atnextpar{\the\toks@}%
 \endgroup
}

\def\ateverypar#1{%
 \evp@setup
 \begingroup
  \expandafter\toks@\expandafter{\@ateverypar #1}%
  \xdef\lib@gtmp{\the\toks@}%
 \endgroup
 \let\@ateverypar\lib@gtmp
 \let\lib@gtmp\undefined
}

\restoreatcode

\saveatcode \makeatletter

\newbox\risbox
\def\setrisbox{\setbox\risbox}
\def\riswidth{\wd\risbox}

\def\rismove(#1,#2){%
 \setrisbox\vbox{%
  \dimen@=#2\relax
  \vskip\dimen@
  \hbox{\dimen@=#1\relax\hskip\dimen@\box\risbox\hskip-\dimen@}%
  \vskip-\dimen@
 }%
}

\let\global@ris@options\lib@empty

\def\resetrisoptions{\let\global@ris@options\lib@empty}

\def\addrisoption#1{%
  \expandafter\toks@\expandafter{\global@ris@options #1\relax}%
  \edef\global@ris@options{\the\toks@}%
}

\def\ris@options{} 

\def\atnextris#1{%
  \expandafter\toks@\expandafter{\ris@options #1\relax}%
  \edef\ris@options{\the\toks@}%
}

\newif\ifrisleft

\def\rischeckpage{%
  \dimen@\ht\risbox
  \advance\dimen@\dp\risbox
  \ifdim\dimen@<\risminheight \dimen@\risminheight \fi
  \advance\dimen@\pagetotal
  \advance\dimen@\pagedepth
  \advance\dimen@ -\pageshrink
  \ifdim\dimen@>\pagegoal
    \ifrisleft \risleftfalse \else \rislefttrue \fi
  \fi
}

\def\ris{\begingroup
  \risleftfalse
  \tempdimen\risspace \risspace 1em
  \temptoks\risadjust \risadjust{}%
  \temptoks\risadjustabove \risadjustabove{}%
  \temptoks\risadjustbelow \risadjustbelow{}%
  \tempdimen\risminheight \risminheight\z@
  \tempcount\risshowbox \risshowbox\z@
  \global@ris@options
  \ris@options
  \futurelet\lib@next\ris@checkopt
}

\def\ris@parseopt{%
  \futurelet\lib@next\ris@checkopt
}

\def\ris@checkopt{%
 \expandafter\if\space\noexpand\lib@next
   \expandafter\ris@checkopt@space
  \else
   \expandafter\ris@checkopt@
  \fi
}
\def\ris@checkopt@space{\afterassignment\ris@parseopt\let\lib@next= }

\def\ris@checkopt@{%
 \ifx\lib@next[\expandafter\ris@getopt\else\expandafter\ris@\fi
}

\def\ris@getopt[#1]{#1\ris@}

\def\ris@{%
  \ris@breakpar{\the\risadjustabove\relax\the\risadjust\relax}%
  \risadjust{}
  \tempdimen\ris@totalspace
  \ris@totalspace\riswidth
  \advance\ris@totalspace\risspace
  \tempdimen\ris@hsize
  \ris@hsize\hsize
  \advance\ris@hsize-\ris@totalspace
  \ifrisleft
    \global\parshape 1 \ris@totalspace \ris@hsize
    \ateverypar{\global\parshape 1 \ris@totalspace \ris@hsize}%
  \else
    \global\parshape 1 \z@ \ris@hsize
    \ateverypar{\global\parshape 1 \z@ \ris@hsize}%
  \fi
  \ifhmode \ris@insert \else \atnextpar{\ris@insert}\fi
  \no@page@breaks
  \ignorespaces
}

\def\ris@breakpar#1{%
  \ifhmode{%
   \parfillskip\z@\parskip\z@\tolerance\@m\par #1\parindent\z@\leavevmode
  }\else #1\fi
}

\def\no@page@breaks{%
  \ifx\tex@par@npb\undefined
    \let\tex@par@npb\par
    \let\tex@penalty@npb\penalty
    \interlinepenalty\@M
    \predisplaypenalty\@M
    \postdisplaypenalty\@M
    \def\par{\tex@par@npb\tex@penalty@npb\@M\relax}%
    \def\penalty{\ifvmode\count@\else\tex@penalty@npb\fi}%
  \fi
}

\def\end@no@page@breaks{%
  \ifx\tex@par@npb\undefined \else
    \let\par\tex@par@npb
    \let\penalty\tex@penalty@npb
    \let\tex@par@npb\undefined
    \let\tex@penalty@npb\undefined
  \fi
}

\def\ris@insert{{%
 \ifx\ris@totalspace\undefined
  \errmessage{Side insertion with no text}
 \else
  \ifnum\risshowbox>0
    \setrisbox\vbox{%
      \vskip-.4pt\hrule 
      \hbox{\hskip-.4pt\vrule \box\risbox \vrule\hskip-.4pt}%
      \hrule\vskip-.4pt
    }%
  \fi
  \setbox\z@\lastbox
  \hbox to\z@{%
   \ifrisleft \kern-\ris@totalspace
   \else \kern\ris@hsize \kern\risspace \fi
   \vbox to\z@{\kern -1ex \box\risbox\vss}\hss
  }%
  \box\z@
 \fi
}}

\def\endris{\end@no@page@breaks\futurelet\lib@next\endris@}

\def\endris@{%
 \ifx\lib@next\par
   \par \the\risadjustbelow\relax\the\risadjust\relax
 \else
   \ris@breakpar{\the\risadjustbelow\relax\the\risadjust\relax}%
 \fi
 \endgroup
 \global\parshape\z@
 \let\ris@options\lib@empty 
}

\def\epsrisbox#1{\hbox{%
  \testfile{#1}\ifyes
    \def\lib@tmp{\epsfbox{#1}}%
  \else
    \message{^^JWarning: file #1 not found^^J}%
    \global\epsfxsize\z@
    \global\epsfysize\z@
    \let\lib@tmp\relax
  \fi
  \lib@tmp
}}

\restoreatcode

\addrisoption{%
  \risminheight 3\baselineskip
}

\addrisoption{%
  \pretolerance=100 \tolerance=2000
  \hyphenpenalty=10 \exhyphenpenalty=10
}

\input cells.sty

\usepackage[text={170mm,250mm},centering]{geometry}

\theoremstyle{plain}
\newtheorem{thm}{\indent\normalfont{T\,h\,e\,o\,r\,e\,m}\,}

\def\lemma{\ifvmode{\smallbreak\par}\fi L\,e\,m\,m\,a}

\begin{document}

\title{Points on a line, shoelace and dominoes}
\author{A.\,Khrabrov \and K.\,Kokhas} 
\date{St.-Petersburg University, Russia}

\maketitle

\begin{abstract}
In this survey we consider numerous known and unknown combinatorial realizations of the sequence 
A079487 and basic facts about it.
\end{abstract}

The following problem of N.~Filonov was proposed at St.~Petersburg mathematical olympiad in 2014.

\smallskip
\begin{narrower}\noindent\slshape
40 points are marked on each of two parallel lines.
They are split onto 40 pairs such that segments that join points in a pair do not intersect.
(In particular, no endpoint of a segment lies on any other segment.)
Prove that the number of these matchings is less than $3^{39}$. \par
\end{narrower}

\noindent
A sequence that describes this number of matchings has a lot of combinatorial realizations.
Their diversity is amazing.
Of course, the number of its realizations is not as big (yet?) as for Catalan numbers,
but the sequence itself is a couple of hundreds years younger than Catalan sequence.
Below  we give a survey of numerous known and unknown combinatorial realizations of this sequence. 


\paragraph{Triangle $\text{\bfseries\itshape a}_{\text{\bfseries\itshape {k,n}}}$.}
Consider the following combinatorial construction.
Given two parallel lines, $k$ points are marked on the first line, and $n$ points are marked on the second line.
The points are split into pairs such that segments that join points in a pair do not intersect.
(In particular, no endpoint of a segment lies on any other segment.
The picture obtained by this construction we will call a \emph{configuration} (of points and segments)
or a \emph{partition} (of points into pairs). Denote the number of partition by $a_{k,n}$.
For example $a_{2,4}=4$, as we can see on the following picture. 

\begin{center}
\setlength{\unitlength}{1.5pt}
\begin{picture}(180,12)
\matrixput(0,10)(50,0){4}(10,0){2}{\circle*{2}}
\matrixput(0,0)(50,0){4}(10,0){4}{\circle*{2}}
\path(1,0)(9,0)  \path(21,0)(29,0)  \path(1,10)(9,10)  
\path(50,1)(50,9)  \path(71,0)(79,0)  \path(60,1)(60,9)
\path(100,1)(100,9)  \path(110,0)(120,0)  \path(130,0)(110,10)
\path(150,0)(160,0)  \path(170,0)(150,10) \path(180,0)(160,10)
\end{picture}
\end{center}

\noindent
Assume that $a_{k,n}=0$, if $k+n$ is odd or at least one of numbers $k$, $n$ is negative.
It is convenient to arrange the nonzero numbers $a_{k,n}$ in a triangle-shaped array:
\begin{equation}
{\footnotesize\setcounter{MaxMatrixCols}{20}
\begin{matrix}
 & & &  & 1         \\
 & & & 1& 1& 1      \\
 & &1& 2& 2& 2&1    \\
 &1&3& 4& 5& 4&3&1  \\
1&4&7&10&11&10&7&4&1\\
\end{matrix}\qquad
\begin{matrix}\arraycolsep=0pt
 & &        &        & a_{0,0}         \\
 & &        & a_{2,0}& a_{1,1}& a_{0,2}\\
 & & a_{4,0}& a_{3,1}& a_{2,2}& a_{1,3}&a_{0,4}    \\
 & &        &        &\dots\\
 & &        &        &\dots
\end{matrix}}
\label{egn:main_triangle}
\end{equation}

\begin{figure}
\kern-2mm
\begin{align*}
&
\hbox{\setlength{\unitlength}{1.5pt}\footnotesize
\begin{picture}(56,30)(-5,-20)
\multiput(0,-5)(10,0){5}{\circle*{2}}
\multiput(0,5)(10,0){3}{\circle*{2}}
\put(0,5){\line(1,0){10}}
\put(-5,9){$A_1$}\put(-5,-12){$B_1$}\put(19,9){$\dots$}
\put(6,9){$A_2$}\put(6,-12){$B_2$}
\put(24,-12){$\dots$}\put(38,-12){$B_{n+1}$}
\put(16,-21){\normalsize $a_{k-1,n+1}$}
\end{picture}}
\quad
&&
\hbox{\setlength{\unitlength}{1.5pt}\footnotesize
\begin{picture}(55,30)(-5,-20)
\multiput(0,-5)(10,0){5}{\circle*{2}}
\multiput(0,5)(10,0){3}{\circle*{2}}
\put(0,-5){\line(0,1){10}}
\put(-5,9){$A_1$}\put(-5,-12){$B_1$}\put(19,9){$\dots$}
\put(6,9){$A_2$}\put(6,-12){$B_2$}
\put(24,-12){$\dots$}\put(38,-12){$B_{n+1}$}
\put(20,-21){\normalsize $a_{k,n}$}
\end{picture}}
\quad
&&
\hbox{\setlength{\unitlength}{1.5pt}\footnotesize
\begin{picture}(55,25)(-5,-20)
\multiput(0,-5)(10,0){5}{\circle*{2}}
\multiput(0,5)(10,0){3}{\circle*{2}}
\put(0,-5){\line(1,0){10}}
\put(-5,9){$A_1$}\put(-5,-12){$B_1$}
\put(6,9){$A_2$}\put(6,-12){$B_2$}\put(19,9){$\dots$}
\put(24,-12){$\dots$}\put(38,-12){$B_{n+1}$}
\put(16,-21){\normalsize $a_{k+1,n-1}$}
\end{picture}}
\\[5pt]
&&&
\hbox{\setlength{\unitlength}{1.5pt}\footnotesize
\begin{picture}(50,36)(-5,-20)
\multiput(0,-5)(10,0){5}{\circle*{2}}
\multiput(0,5)(10,0){3}{\circle*{2}}
\put(0,-5){\line(1,0){10}}
\put(0,5){\line(1,0){10}}
\put(-5,9){$A_1$}\put(-5,-12){$B_1$}\put(19,9){$\dots$}
\put(6,9){$A_2$}\put(6,-12){$B_2$}
\put(24,-12){$\dots$}\put(38,-12){$B_{n+1}$}
\put(16,-21){\normalsize $a_{k-1,n-1}$}
\end{picture}}
\end{align*}

\centering
\caption{Number of partitions when we remove one segment.} 
\label{ris:rek_sootn}
\end{figure}
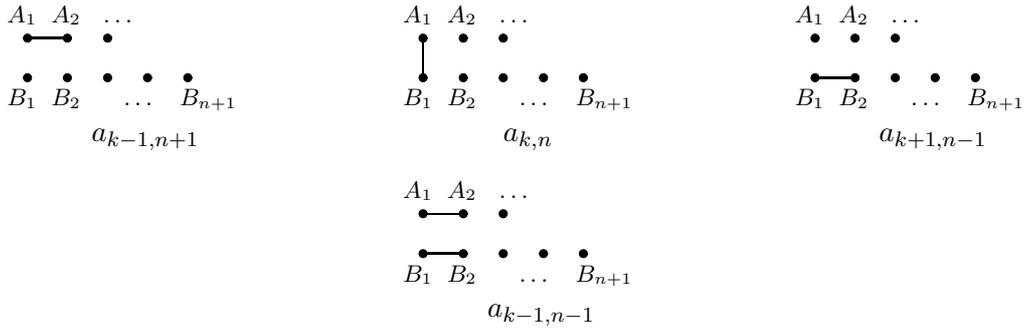

\paragraph{Recurrence reltions.} 
Let us express the number $a_{k+1,n+1}$ via the numbers $a_{i,j}$ with a smaller sum of indices.
Denote points on the upper line by $A_1$, $A_2$,~\dots, $A_{k+1}$  (from left to right)
and on the lower line by $B_1$, $B_2$,~\dots, $B_{n+1}$.
The number of partitions containing the edge $A_1B_1$ equals $a_{k,n}$, 
the number of partitions containing the edge $A_1A_2$ equals $a_{k-1,n+1}$, 
and the number of partitions containing the edge $B_1B_2$ equals $a_{k+1,n-1}$. 
The sum of these three numbers  is equal to $a_{k+1,n+1}+a_{k-1,n-1}$ because
we count twice the partitions containing both segments $A_1A_2$ and $B_1B_2$
(fig.~\ref{ris:rek_sootn}). 
Therefore we have the following remarkable recurrence
\begin{equation}
a_{k+1,n+1}= a_{k,n} + a_{k-1,n+1} + a_{k+1,n-1} - a_{k-1,n-1}.
\label{eqn:4term-recurrence}
\end{equation}
For example the number 11 on the bottom side of the triangle \eqref{egn:main_triangle} equals the sum $4+5+4-2$.
Since we assume the numbers to be equal 0 for negative indices, the recurrence~\eqref{eqn:4term-recurrence} remains valid
on the boundary of the triangle. It does not work in the case $k=n=-1$ only, i.e. for the topmost 1 of the triangle.

Now we will prove one more recurrence relation for $a_{k,n}$:
\begin{equation}
a_{k,n}=a_{k-2,n} + (a_{k-1,n-1}+a_{k-1,n-3}+a_{k-1,n-5}+\ldots).
\label{eqn:long-recurrence}
\end{equation}
The term $a_{k-2,n}$ is equal to the number of configurations containing the segment  $A_1A_2$.
All other configurations contain a segment that joins $A_1$ with some point on the bottom line;
the sum in parentheses counts these configurations.
Indeed, the point $A_1$ can be joined with points $B_3$, $B_5$, \dots \ only,
because otherwise there are odd number of points below the line $A_1B_k$.
If we join $A_1$ and $B_{2i+1}$ ($i=0$, 1, 2 \dots), then the configuration must contain
also the segments $B_1B_2$, \dots, $B_{2i-1}B_{2i}$, hence the number of configurations in this case equals $a_{k-1,n-2i-1}$.

Due to the recurrence \eqref{eqn:long-recurrence} and the symmetry $a_{k,n}=a_{n,k}$
it is not difficult to check that numbers $a_{k,n}$ increase if we move along the row towards its center.

\setrisbox\vbox{\hsize=52pt\noindent
\begin{picture}(54,24)
\multiput(0,-2)(10,0){6}{\line(0,1){24}}
\multiput(-2,0)(0,10){3}{\line(1,0){54}}
\Thicklines
\put(0,10){\vector(1,0){10}}
\put(20,10){\vector(1,1){10}}
\put(40,10){\vector(1,-1){10}}
\end{picture}}

\paragraph{Motzkin peakless paths.} 

\ris
\emph{Motzkin path} is an (oriented) lattice path with steps depicted on the right figure.
We mean that the beginning of this path is in the point $(0,0)$. \emph{Peakless} Motzkin path is
a path that has no peaks, i.e. the fragment of the form
\endris
\lower3pt\hbox{\vbox{\hsize=23pt\noindent 
\begin{picture}(24,14)(-2,-2)
\linethickness{0.25pt}
\multiput(0,-2)(10,0){3}{\line(0,1){14}}
\multiput(-2,0)(0,10){2}{\line(1,0){24}}
\Thicklines
\put(0,0){\vector(1,1){10}}
\put(10,10){\vector(1,-1){10}}
\end{picture}}} . 
Denote by $m_{k,n}$ the number of peakless Motzkin paths from the point $(0,0)$ to the point $(k,n)$. 
It is clear that $-k\leq n \leq k$.
On the next figure you can see all peakless Motzkin paths from the point $(0,0)$ to the point $(3,1)$,
so $m_{3,1}=4$.

\begin{center}
\begin{picture}(200,24)(-3,-2)
\linethickness{0.25pt}
\matrixput(0,-2)(10,0){4}(55,0){4}{\line(0,1){24}}
\matrixput(-2,0)(0,10){3}(55,0){4}{\line(1,0){34}}
\Thicklines
\path(0,10)(10,0)(30,20)
\path(55,10)(75,10)(85,20)
\path(110,10)(120,20)(140,20)
\path(165,10)(175,10)(185,20)(195,20)
\end{picture}
\end{center}



\begin{thm}
$m_{k,n}=a_{k-n,k+n}$. 
\end{thm}

\proof 
What is the last step of a Motzkin path that comes to the point $(k,n)$? 
It is clear that can it start from the three points only: $(k-1,n-1)$, or $(k-1,n)$, or $(k-1,n+1)$.
In the last case we should keep in mind that since the peak in the point $(k-1,n+1)$ is prohibited,
the previous step can not start from the point $(k-2,n)$.
Thus we obtain a recurrence
$$
m_{k,n}= m_{k-1,n-1} + m_{k-1,n} + m_{k-1,n+1}- m_{k-2,n}.
$$
This is a recurrence \eqref{eqn:4term-recurrence} up to a change of variables.
\endproof

\paragraph{A sum of binomial coefficients.} 

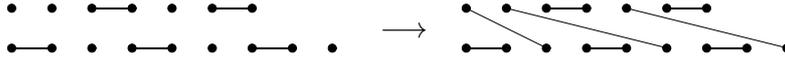
\begin{figure}
\begin{center}
\setlength{\unitlength}{1.5pt}
\noindent
\begin{picture}(75,13)
\multiput(0,0)(10,0){9}{\circle*{2}}
\multiput(0,10)(10,0){7}{\circle*{2}}
\put(0,0){\line(1,0){10}}
\put(30,0){\line(1,0){10}}
\put(60,0){\line(1,0){10}}
\put(20,10){\line(1,0){10}}
\put(50,10){\line(1,0){10}}
\end{picture}
\qquad\lower-4pt\hbox{$\longrightarrow$}\quad
\setlength{\unitlength}{1.5pt}
\begin{picture}(80,13)
\multiput(0,0)(10,0){9}{\circle*{2}}
\multiput(0,10)(10,0){7}{\circle*{2}}
\put(0,0){\line(1,0){10}}
\put(30,0){\line(1,0){10}}
\put(60,0){\line(1,0){10}}
\put(20,10){\line(1,0){10}}
\put(50,10){\line(1,0){10}}
\put(0,10){\line(2,-1){20}}
\put(10,10){\line(4,-1){40}}
\put(40,10){\line(4,-1){40}}
\end{picture}
\end{center}

\kern-3mm
\centering
\caption{
Horizontlal segments determine a partition onto pairs uniquely.
}
\label{ris:horrebra}
\end{figure}

Let  numbers $k$ and $n$ be of the same parity.
Consider an arbitrary configuration of $k$ on the upper line and $n$ points on the lower line.
Then remove all the segments whose endpoints belong to the different lines. 
The remaining set of horizontal segments and disjoined points determines uniquely 
the initial partition of points onto pairs (see fig.~\ref{ris:horrebra}).
Indeed, the leftmost point on the upper line must be joined with the leftmost point on the lower line,
the second (from the left) point must be joined with the second (from the left) point on the lower line, etc.,
because otherwise the segments will intersect.


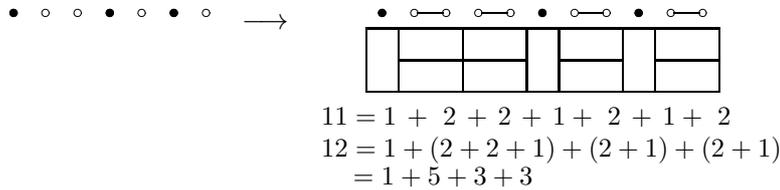
\begin{figure}[h]
\begin{center}
\setlength{\unitlength}{1.2pt}
\noindent
\begin{picture}(50,35)(0,-30)
\multiput(0,0)(10,0){7}{\circle{2.25}}
\put(0,0){\circle*{2.25}}
\put(30,0){\circle*{2.25}}
\put(50,0){\circle*{2.25}}
\end{picture}
\qquad
\raise30pt\hbox{$\longrightarrow$}
\ \
\begin{picture}(132,5)(-20,-30)
\multiput(0,0)(10,0){11}{\circle{2.25}}
\put(0,0){\circle*{2.25}}
\put(50,0){\circle*{2.25}}
\put(80,0){\circle*{2.25}}
\put(11,0){\line(1,0){8}}
\put(31,0){\line(1,0){8}}
\put(61,0){\line(1,0){8}}
\put(91,0){\line(1,0){8}}
\put(-5,-25){\box1}
\put(-19,-35){\small $11=1\,+\,\,2\,+\,2\,+\,1+\,\,2\,+\,1+\,\,2$}
\put(-19,-45){\small $12=1+(2+2+1)+(2+1)+(2+1)$}
\put(-9,-54){\small $=1+5+3+3$}
\end{picture}
\end{center}

\kern4mm
\centering
\caption{\vtop{\hsize=330pt A colouring in two colors determines uniquely horizontal segments, 
partitions onto dominoes and compositions.}}
\label{ris:conf-domino-compos}

\end{figure}

Consider an auxiliary problem.
We call a $(j,\ell)$-configuration on a line a set of $\ell$ non-intersecting segments and 
$j$ points that do not lay on these segments. We color all endpoints of segments in white 
and all other alone points in black.
Let us find the number of $(j,\ell)$-configuration on a line.
To do this we will depict configurations on the following order.
First, mark $j+\ell$ points on the line, choose $\ell$ of the points and color them in white.
Then, for every white point add a new white point in its small left neighborhood and join these points by a segment
(see fig.~\ref{ris:conf-domino-compos}).
It is clear that different choices of $\ell$ white points give us different configurations.
Therefore, the number of configurations equals the number of choices, i.e. $\binom{j+\ell}\ell$.

\setbox1\hbox{\cellsize=12pt\cells{
 _ _ _ _ _ _ _ _ _ _ _
| |_ _|_ _| |_ _| |_ _|
|_|_ _|_ _|_|_ _|_|_ _|}}

Let $k$ be the total number of points in a $(j,\ell)$-configuration, i.e. $k=j+2\ell$.
Then the last binomial coefficient can be written in the form 
$\binom{(k+j)/2}{(k-j)/2}$ or $\binom{(k+j)/2}j$. 
Returning to the initial problem about points on two lines we see that
the number of ways to choose several horizontal segments on the upper line 
($k$ points and $j$ of them are ``alone'') equals $\binom{(k+j)/2}j$,
and the similar number for the lower line ($n$ points, $j$ of them are alone) equals $\binom{(n+j)/2}j$,
here $k$ and $n$ are of the same parity due problem conditions, and $j$ must have the same parity, too, 
because the non-alone points are split onto pairs.
Thus we obtain an explicite formula for $a_{k,n}$:
$$
a_{k,n} = \sum_{\substack{0\leq j\leq\min\{k,n\} \\  j\equiv k\pmod{2}}} 
C_{\frac{k+j}{2}}^{\frac{k-j}{2}}C_{\frac{n+j}{2}}^{\frac{n-j}{2}} =
\sum_{\substack{0\leq j\leq\min\{k,n\} \\  j\equiv k\pmod{2}}} 
C_{\frac{k+j}{2}}^{j}C_{\frac{n+j}{2}}^{j}.
$$
For $k=n$ the formula looks simpler, because 
$\frac{n-j}{2}=\frac{k-j}{2}=\ell$ and  $\frac{n+j}{2}=\frac{k+j}{2}=k-\ell$ in this case, 
and so
$$
a_{n,n} = \sum_{\ell=0}^{[n/2]} \big(C_{n-\ell}^\ell\big)^2 .
$$

\paragraph{Domino tilings.}
Consider a domino tilings of rectangles $2\times k$ and $2\times n$.
Let $d_{k,n}$ be a number of those tilings that have equal numbers of vertical domino. 
For example the following table contains domino tilings of rectangles $2\times 3$ and $2\times 5$
with 1 and 3 vertical dominoes. From this table we see that  $d_{3,5}=6+4=10$.

\bigskip
\begin{center}
{\footnotesize\extrarowheight=8pt\noindent\ \
\begin{tabular}{|c|c|c|c|}
\hline
$j$ &  $2\times 3$ & $2\times 5$ & number of variants\\  \hline
1 & \cells{
 _ _ _
| |_ _|
|_|_ _|
}\,, \  \cells{
 _ _ _
|_ _| |
|_ _|_|
} &  \cells{
 _ _ _ _ _
| |_ _|_ _|
|_|_ _|_ _|
}\,, \  \cells{
 _ _ _ _ _
|_ _| |_ _|
|_ _|_|_ _|
}\,, \  \cells{
 _ _ _ _ _
|_ _|_ _| |
|_ _|_ _|_|
} &  $2\cdot3=6$ \\  \hline
3 & \cells{
 _ _ _
| | | |
|_|_|_|
} & 
\cells{
 _ _ _ _ _
| | | |_ _|
|_|_|_|_ _|
}\,, \cells{
 _ _ _ _ _
| | |_ _| |
|_|_|_ _|_|
}\,, \cells{
 _ _ _ _ _
| |_ _| | |
|_|_ _|_|_|
}\,, \cells{
 _ _ _ _ _
|_ _| | | |
|_ _|_|_|_|
} & $1\cdot4=4$ \\ \hline
\end{tabular}}
\end{center}

\bigskip
\begin{thm}
$d_{k,n}=a_{k,n}$. 
\label{thm:a=dominoshki}
\end{thm}

\proof Has written in the previous paragraph. 
We have to replace the words ``$(j,\ell)$-configu\-ra\-tion'' by
``domino tiling of the rectangle $2\times(j+2\ell)$'';
the words ``horizontal segment'' by words ``pair of horizontal dominoes''
and the words ``black point'' by ``vertical domino'', see fig.~\ref{ris:conf-domino-compos}.
\endproof

\paragraph{Fibonacci numbers.}

Let $f_m$ be the sum of the numbers in $m$-th row of our triangle~\eqref{egn:main_triangle}, i.e.
$$
f_m = a_{2m,0}+a_{2m-1,1}+a_{2m-2,2}+\ldots+a_{1,2m-1}+a_{0,2m} .
$$
Let us sum up recurrences~\eqref{eqn:4term-recurrence} over all numbers $k$ and $n$ with $k_n=2m+2$.
In the l.h.s. we obtain a sum of all elements of $(m+1)$-th row.
In the r.h.s we obtain triple sum of all elements of $m$-th row (with plus sign) and
sum of all elements of $(m-1)$-th row (with minus sign).
Hence  $f_{m+1} = 3f_m-f_{m-1}$. But this recurrence is also valid for Fibonacci numbers $F_{2m}$:
$$
F_{2m+2} = F_{2m+1}+F_{2m} = 2F_{2m}+F_{2m-1} =
2F_{2m}+(F_{2m}-F_{2m-2}) = 3F_{2m} - F_{2m-2} .
$$
Since the sums of the first and the second rows of the triangle equal $F_2=1$ and $F_4=3$, 
we have the equality $f_m=F_{2m}$. So the sum of elements of $m$-th row in the triangle 
equals $2m$-th Fibonacci number.

From this observation we immediately obtain the estimation $a_{k,n}\leq F_{n+k}$.
Using Binet formula for Fibonacci numbers we can rewrite it in a more concrete form.
By Binet formula
$$
F_m=\frac{\Bigl(\frac{1+\sqrt{5}}2\Bigr)^m-\Bigl(\frac{1-\sqrt{5}}2\Bigr)^m}{\sqrt{5}}.
$$
The second term in the numerator has small absolute value and negative for even $n$.
By omitting this term we slightly increse the r.h.s. and obtain
\begin{equation}
a_{k,n}< \frac1{\sqrt{5}} \biggl(\frac{1+\sqrt{5}}2\biggr)^{n+k}\,.
\label{eqn:ocenka_a_nk_sverhu}
\end{equation}

\paragraph{Fences.}

Consider the following bipartite oriented graph $Z_{2n}$ whose parts contain $n$ vertices.
We call this graph \emph{a fence}.

\begin{center}
\begin{picture}(100,20)
\multiput(0,0)(20,0){5}{\circle{3}}
\multiput(10,20)(20,0){5}{\circle{3}}
\multiput(9,18)(20,0){5}{\vector(-1,-2){8}}
\multiput(11,18)(20,0){4}{\vector(1,-2){8}}
\end{picture}
\end{center}

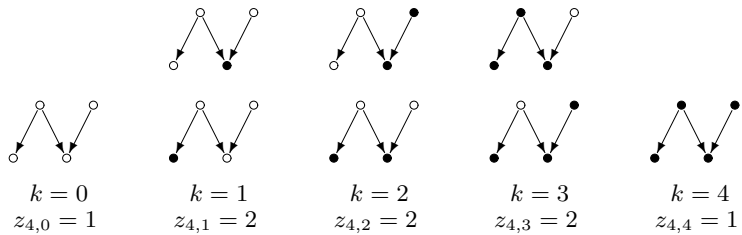
\begin{figure}[b]
\begin{center}\footnotesize
\setlength{\unitlength}{1pt}
\begin{picture}(280,81)(0,-26)
\multiput(0,0)(20,0){2}{\circle{3}}
\multiput(10,20)(20,0){2}{\circle{3}}
\multiput(9,18)(20,0){2}{\vector(-1,-2){8}}
\put(11,18){\vector(1,-2){8}}
\put(6,-16){$k=0$}
\put(0,-26){$z_{4,0}=1$}
\put(60,0){
\put(0,0){\circle*{3}}
\put(20,0){\circle{3}}
\multiput(10,20)(20,0){2}{\circle{3}}
\multiput(9,18)(20,0){2}{\vector(-1,-2){8}}
\put(11,18){\vector(1,-2){8}}
\put(6,-16){$k=1$}
\put(0,-26){$z_{4,1}=2$}}
\put(60,35){
\put(0,0){\circle{3}}
\put(20,0){\circle*{3}}
\multiput(10,20)(20,0){2}{\circle{3}}
\multiput(9,18)(20,0){2}{\vector(-1,-2){8}}
\put(11,18){\vector(1,-2){8}}
}
\put(120,0){
\put(0,0){\circle*{3}}
\put(20,0){\circle*{3}}
\multiput(10,20)(20,0){2}{\circle{3}}
\multiput(9,18)(20,0){2}{\vector(-1,-2){8}}
\put(11,18){\vector(1,-2){8}}
\put(6,-16){$k=2$}
\put(0,-26){$z_{4,2}=2$}}
\put(120,35){
\put(0,0){\circle{3}}
\put(20,0){\circle*{3}}
\put(10,20){\circle{3}}
\put(30,20){\circle*{3}}
\multiput(9,18)(20,0){2}{\vector(-1,-2){8}}
\put(11,18){\vector(1,-2){8}}
}
\put(180,0){
\put(0,0){\circle*{3}}
\put(20,0){\circle*{3}}
\put(10,20){\circle{3}}
\put(30,20){\circle*{3}}
\multiput(9,18)(20,0){2}{\vector(-1,-2){8}}
\put(11,18){\vector(1,-2){8}}
\put(6,-16){$k=3$}
\put(0,-26){$z_{4,3}=2$}}
\put(180,35){
\put(0,0){\circle*{3}}
\put(20,0){\circle*{3}}
\put(10,20){\circle*{3}}
\put(30,20){\circle{3}}
\multiput(9,18)(20,0){2}{\vector(-1,-2){8}}
\put(11,18){\vector(1,-2){8}}
}
\put(240,0){
\multiput(0,0)(20,0){2}{\circle*{3}}
\multiput(10,20)(20,0){2}{\circle*{3}}
\multiput(9,18)(20,0){2}{\vector(-1,-2){8}}
\put(11,18){\vector(1,-2){8}}
\put(6,-16){$k=4$}
\put(0,-26){$z_{4,4}=1$}}
\end{picture}
\end{center}

\caption{Closed sets of the graph $Z_4$.}
\label{fig:zamkn-nabory-Z4}
\end{figure}

The set $A$ of vertices of this graph is called \emph{closed}, if we can not leave this set moving along arrows.
In other words it satisfies the property: if  $x\to y$ is an edge of the graph and $x\in A$, then $y\in A$.
Denote by $z_{2n,k}$ the number of $k$-element closed sets of vertices in the graph  $Z_{2n}$ ($0\leq k\leq 2n$).
In fig.~\ref{fig:zamkn-nabory-Z4} the vertices of closed subsets of the graph  $Z_4$ are colored in black.

\pagebreak

\begin{thm}
$z_{2n,k}=a_{2n-k,k}$. 
\end{thm}

\proof See fig.~\ref{ris:bij1}.
\endproof

\begin{figure}[h]
\begin{center}\footnotesize
\setlength{\unitlength}{1pt}
%
%
\begin{picture}(150,20)
\multiput(0,0)(20,0){2}{\circle{3}}
\multiput(40,0)(20,0){4}{\circle*{3}}
\put(120,0){\circle{3}}
\put(140,0){\circle*{3}}
\multiput(10,20)(20,0){2}{\circle{3}}
\put(50,20){\circle*{3}}\put(90,20){\circle*{3}}\put(150,20){\circle*{3}}
\multiput(110,20)(20,0){2}{\circle{3}}
\put(70,20){\circle{3}}
\multiput(9,18)(20,0){8}{\vector(-1,-2){8}}
\multiput(11,18)(20,0){7}{\vector(1,-2){8}}

\end{picture}

\bigskip
\leavevmode
Mark the vertices that belong to a closed set with black color.

%
%
\begin{picture}(170,41)
\multiput(0,20)(10,0){4}{\circle{3}}
\multiput(2,20)(10,0){3}{\line(1,0){6}}
\multiput(45,0)(10,0){3}{\circle*{3}}
\multiput(47,0)(10,0){2}{\line(1,0){6}}
\put(80,20){\circle{3}}
\multiput(95,0)(10,0){3}{\circle*{3}}
\multiput(97,0)(10,0){2}{\line(1,0){6}}
\multiput(130,20)(10,0){3}{\circle{3}}
\multiput(132,20)(10,0){2}{\line(1,0){6}}
\multiput(165,0)(10,0){2}{\circle*{3}}
\put(167,0){\line(1,0){6}}
\end{picture}

\bigskip
\leavevmode
\vbox{\hsize=450pt\noindent 
Arrange the vertices onto two parallel lines.
On the upper line the leftmost (possibly empty) group of vertices is even,
all other groups are odd;  on the lower line the rightmost group is even, all other groups are odd.
}

\bigskip
%
%
\begin{picture}(170,31)
\multiput(0,20)(10,0){4}{\circle{3}}
\multiput(2,20)(20,0){2}{\line(1,0){6}}
\multiput(45,0)(10,0){3}{\circle*{3}}
\put(47,0){\line(1,0){6}}
\path(68,3)(78,17)
\put(80,20){\circle{3}}
\multiput(95,0)(10,0){3}{\circle*{3}}
\put(97,0){\line(1,0){6}}
\path(118,3)(128,17)
\multiput(130,20)(10,0){3}{\circle{3}}
\put(142,20){\line(1,0){6}}
\multiput(165,0)(10,0){2}{\circle*{3}}
\put(167,0){\line(1,0){6}}
\end{picture}

\bigskip
\leavevmode
\vbox{\hsize=330pt\noindent 
Join the last vertex in every odd group on lower line with the first vertex of the next group on the upper line.
All other vertices in groups split onto pairs.
}
\end{center}


\centering
\caption{\vtop{\hsize=270pt A bijection between the closed sets in graph 
$Z_{2n}$ and partitions of points onto pairs.
}}
\label{ris:bij1}
\end{figure}

\paragraph{Fibonacci numbers with odd indices.} 

An inspective reader probably has been concerned for some time about the issue:
why the row sums of the triangle  \eqref{egn:main_triangle} are given by Fibonacci numbers with even indices only?
And where we can find in this topic Fibonacci numbers with odd indices?
Here you are.
Consider an odd fence $Z_{2n+1}$ containing $n$ vertices in the upper part, and $n+1$ vertices in the lower part.

\begin{equation}
\begin{picture}(122,20)
\multiput(0,0)(20,0){6}{\circle{3}}
\multiput(10,20)(20,0){5}{\circle{3}}
\multiput(9,18)(20,0){5}{\vector(-1,-2){8}}
\multiput(11,18)(20,0){5}{\vector(1,-2){8}}
\put(104,-3){\small $A$}
\put(93,19){\small $B$}
\end{picture}
\label{eqn:ris-zabora}
\end{equation}

\noindent
The definition of the closed set for this fence is the same as for even fence.
Denote by $z_{2n+1,k}$ the number of $k$-element closed sets of vertices in the graph  $Z_{2n+1}$ ($0\leq k\leq 2n+1$).
Remark that the numbers $z_{2n+1,k}$ generally speaking are not symmetric, i.e. usually
$z_{2n+1,k}\ne z_{2n+1,2n+1-k}$.

It is easy to check that the following recurrence relations hold:
\begin{align*}
z_{2n, k} &= z_{2n-1,k} + z_{2n-2,k-2},
\\
z_{2n+1, k} &= z_{2n,k-1} + z_{2n-1,k}.
\end{align*}
Let us check for example the second equality. Let $A$ and $B$ be the rightmost vertices of the fence $Z_{2n+1}$
as in the diagram~\eqref{eqn:ris-zabora}.
Consider an arbitrary closed set with $k$ vertices. If $A$ belongs to this set, 
then remove it and we obtain a $(k-1)$-element closed set of the fence $Z_{2n}$.
If $A$ does not belong to the closed set, then $B$ does not belong to this set, too, 
and we remove both vertices and obtain $k$-element set of the fence $Z_{2n-1}$.
It is clear that both types of operations are bijections between the corresponding collections 
of the closed sets and the formula follows.

Consider the union of sequences $z_{n,k}$ for even and odd $n$.  
We can arrange elements of this sequence in the triangle-shaped array similar to Pascal triangle (fig.~\ref{ris:full-treug}).
The recurrence relations written above are illustrated by equalities $10=7+3$ and $5=4+1$ for the numbers
in frames.

\def\circboxed#1{\begin{picture}(8,8)\put(5,2.5){\circle{11}}\put(3,0){#1}\end{picture}}

\begin{figure}
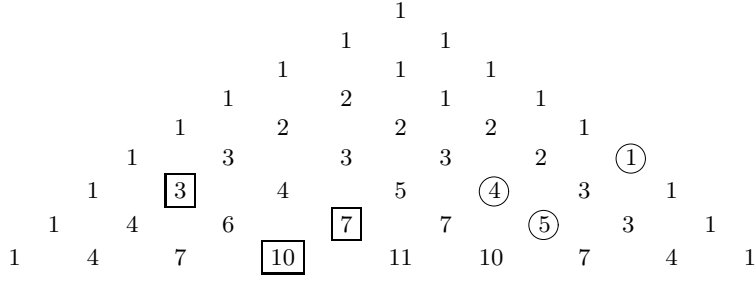

\footnotesize
$$
\setcounter{MaxMatrixCols}{26}
\begin{matrix}
&  &  &  &  &  &  &  &  1                       \\
&  &  &  &  &  &  &  1 &&  1                    \\
&  &  &  &  &  &  1 &&  1 &&  1                 \\
&  &  &  &  &  1 &&  2 &&  1 &&  1              \\
&  &  &  &  1 &&  2 &&  2 &&  2 &&  1           \\
&  &  &  1 &&  3 &&  3 &&  3 &&  2 &&  \circboxed{1}        \\ 
&  &  1 &&  \boxed{3} &&  4 &&  5 &&  \circboxed{4} &&  3 &&  1     \\
&  1 &&  4 &&  6 &&  \boxed{7} &&  7 &&  \circboxed{5} &&  3 &&  1  \\
1 &&  4 &&  7 && \boxed{10} && 11 && 10 &&  7 &&  4 && 1\\
\end{matrix}
$$

\centering
\caption{Triangle $z_{n,k}$}
\label{ris:full-treug}
\end{figure}

Reasoning like that of in the paragraph ``Fibonacci numbers'' shows that  the row sums of this triangle
satisfy the recurrence (and initial conditions) for Fibonacci numbers:  $f_{m+1}=f_{m}+f_{m-1}$.

If we would like to consider odd fences with $n+1$ vertices in the upper part and $n$ vertices in the lower part,
we will obtain a similar recurrences and in fact the same triangle.

\paragraph{Generation function.} 

Consider a formal power series
$$
F(x,y)=\sum_{k,n} a_{k,n}  x^k y^n.
$$
The recurrence \eqref{eqn:4term-recurrence} means that the following 
equality for the function $F(x,y)$ holds
$$ 
(1-x^2 - y^2 +x^2y^2-xy) F(x,y)=1.
$$
Hence 
\begin{equation}
F(x,y)=\frac1{(1-x^2)(1-y^2)-xy}.
\label{eqn-gen-function}
\end{equation}

\paragraph{Compositions.}

Fix a set $S\subset \mathbb{N}$.
\emph{A composition} of a number $n$ is a representation of the number $n$ as a sum in which
each summand belongs to $S$. The representations with different order of summands suppose to be different.

For example a $(j,\ell)$-configuration or a domino tiling of the rectangle  $2\times n$ can be interpreted 
as compositions of $n$ with $S=\{1,2\}$ or as 
compositions of $n+1$ with odd summands (see fig.~\ref{ris:conf-domino-compos}).
This allows us to reformulate theorem \ref{thm:a=dominoshki} and construct two more combinatorial realizations 
of the sequence  $a_{n,k}$ in terms of compositions.
Let $S_1=\{1,2\}$, $S_2=\{1, 3, 5, 7, \ldots\}$ --- odd numbers, $S_3=\mathbb{N}\setminus\{1\}=\{2, 3, 4, \ldots\}$. 

\smallbreak

T\,h\,e\,o\,r\,e\,m \,\ref{thm:a=dominoshki}$'$. The number $a_{k,n}$ equals the number of pairs of 
$S_1$-compositions of numbers $k$ and $n$ with equal number of unities.

T\,h\,e\,o\,r\,e\,m \,\ref{thm:a=dominoshki}$''$. 
The number $a_{k,n}$ equals the number of pairs of $S_2$-compositions of numbers $k+1$ and $n+1$
that have equal number of summands.

\smallskip
Now consider a similar construction, though its relationship to dominoes and configurations is not quite clear.

\smallskip
O\,b\,s\,e\,r\,v\,a\,t\,i\,o\,n. \,The number of $S_1$-compositions of $n$ that contain $\ell$ 2's ($\ell=0$, 1, \dots)
is equal to the number of $S_3$-compositions of $n+2$ with $\ell+1$ summands.
For example there exist 4 $S_1$-compositions of the number 5 that contain one 2
and in the same time there exist 4 $S_3$-compositions of the number 7 that contain two summands:
\begin{gather*}
5=1+1+1+2=1+1+2+1=1+2+1+1=2+1+1+1,
\\
7=2+5=5+2=3+4=4+3.
\end{gather*}
This observation follows from the fact that the number of both compositions equals $\binom{n-\ell}{\ell}$. Prove!

Due to this observation we can change the type of compositions in theorem~\ref{thm:a=dominoshki}$'$.

\smallskip
T\,h\,e\,o\,r\,e\,m \,\ref{thm:a=dominoshki}$'''$. 
Let numbers $n$ and $k$ be of the same parity, $n\geq k$, $s=\frac12(n-k)$.
The number $a_{k,n}$ equals the number of pairs of $S_3$-compositions of numbers $k+2$ and $n+2$
in which the composition of $n+2$ contains $s$ more summands than the composition of $k+2$.


\paragraph{Generation function for compositions} 

Fix a set $S\subset\mathbb{N}$. 
Let $t_n$ be a number of $S$-compositions of $n$, and $T(x)$ be a generation function of the sequence $t_n$.
A very simple formula turns out to hold for this generation function:
\begin{equation}
T(x)=\sum_{n=0}^{+\infty} t_n x^n =    \frac1{1-\sum\limits_{m\in S} x^m}.
\label{eqn:compositions-gen-func}
\end{equation}
Indeed, taking $q=\smash{\sum\limits_{m\in S} x^m}$ and using the formula for geometric series 
$\frac1{1-q}=1+q+q^2+q^3+\ldots$, we obtain
$$
\frac1{1-\sum\limits_{m\in S} x^m} = 1+\Bigl(\sum\limits_{m\in S} x^m\Bigr)+
\Bigl(\sum\limits_{m\in S} x^m\Bigr)^2+\Bigl(\sum\limits_{m\in S} x^m\Bigr)^3+\ldots
$$
Now the reader can expand parentheses and discover a bijection between a set of $S$-compositions of $n$ and 
different occurrences of $x^n$ in this sum, that proves the formula~\eqref{eqn:compositions-gen-func}.

\def\VV{\smash[b]{\lower3.5pt\vbox{\hsize=16pt\noindent 
\begin{picture}(14,14)(-2,-2)
\multiput(0,-2)(10,0){2}{\line(0,1){14}}
\multiput(-2,0)(0,10){2}{\line(1,0){14}}
\Thicklines
\put(0,0){\vector(0,1){10}}
\end{picture}}}}

\def\HH{\lower3.5pt\vbox{\hsize=16pt\noindent 
\begin{picture}(14,14)(-2,-2)
\multiput(0,-2)(10,0){2}{\line(0,1){14}}
\multiput(-2,0)(0,10){2}{\line(1,0){14}}
\Thicklines
\put(0,0){\vector(1,0){10}}
\end{picture}}}

\medskip
E\,x\,a\,m\,p\,l\,e \,1. Consider a function  
$f(x,y)=\frac1{1-x-y}$ and expand it like in the previous reasoning:
$$
\frac1{1-x-y}=1 + (x+y) + (x+y)^2 +(x+y)^3+\ldots
$$
Let us expand parentheses directly without any simplifications, changing the order of summands,  binomial formulae, etc.
Then we obtain terms of the form  $x^ky^m$, and each of them is written as a word with letters $x$ and $y$. 
For example $x^3y^2$ can be obtained as products $xxxyy$, $xyyxx$, etc.
We can interpret such words as ``vector compositions''.
To do this, consider a path on the squared grid which contains a segment \HH for each letter $x$, 
and a segment \VV for each letter $y$. Then the term $x^ky^m$ corresponds to a grid path from the point $(0,0)$ to the point $(k,m)$.
In other words, we represent vector $(k,m)$ as a sums where each summand is a vector $(1,0)$ or $(0,1)$.
These sums are vector compositions. It is clear that the number of these compositions equals $\binom{k+m}k$.
Therefore $f(x,y)$ is a generation function for binomial coefficients:
$$
f(x,y)=\frac1{1-x-y}=\sum_{k,m=0}^{+\infty} C_{m+k}^k x^ky^{m}.
$$

\smallskip
E\,x\,a\,m\,p\,l\,e \,2. Consider a generation function~\eqref{eqn-gen-function}:
$$
F(x,y)=\frac1{1-(x^2+y^2+xy-x^2y^2)}=\sum_k (x^2+y^2+xy-x^2y^2)^k.
$$
Let us ignore for a short time the minus sign before $x^2y^2$.
Then we can expand parentheses and interpret this action as a constructing of vector composition.
Consider an arbitrary monomial $x^ky^m$ in the r.h.s. 
It is equal to the product of multipliers of the form $x^2$, $y^2$, $xy$ or $x^2y^2$.
We put into correspondence to each of these multipliers a directed segment as shown in the picture.
Then the whole monomial is depicted as path from the point $(0,0)$ to $(k,m)$.

\begin{center}
\begin{picture}(24,24)(-2,-2)
\multiput(0,-2)(10,0){3}{\line(0,1){24}}
\multiput(-2,0)(0,10){3}{\line(1,0){24}}
\Thicklines
\put(0,10){\vector(1,0){20}}
\put(7,-12){$x^2$}
\end{picture}
\qquad
\begin{picture}(24,24)(-2,-2)
\multiput(0,-2)(10,0){3}{\line(0,1){24}}
\multiput(-2,0)(0,10){3}{\line(1,0){24}}
\Thicklines
\put(10,0){\vector(0,1){20}}
\put(7,-12){$y^2$}
\end{picture}
\qquad
\begin{picture}(24,24)(-2,-2)
\multiput(0,-2)(10,0){3}{\line(0,1){24}}
\multiput(-2,0)(0,10){3}{\line(1,0){24}}
\Thicklines
\put(0,0){\vector(1,1){10}}
\put(5,-12){$xy$}
\end{picture}
\qquad
\begin{picture}(24,24)(-2,-2)
\multiput(0,-2)(10,0){3}{\line(0,1){24}}
\multiput(-2,0)(0,10){3}{\line(1,0){24}}
\Thicklines
\put(0,0){\vector(1,1){20}}
\put(2,-12){$x^2y^2$}
\end{picture}
\end{center}

\smallskip

\noindent
Thus $a_{k,n}$ equals to the ``number'' of paths from the point $(0,0)$ to $(k,m)$ with specified steps.
We use the word ``number'' in quotes since due to the minus sign we have not pay attention yet
the paths containing $m$ steps of the form  $x^2y^2$ should be equipped with the sign  $(-1)^m$.

In order to cancel these minuses, let us group some paths into pairs.
Consider all the paths that contain a step $x^2y^2$ or two consecutive steps  $x^2$, $y^2$ (in this order).
If the path does not contain a~pair of consecutive steps  $x^2$, $y^2$ before the first occurrence of the step $x^2y^2$ we say that it is of \emph{type A},
otherwise it is of \emph{type~B}.
If we replace in an arbitrary path of type A the first occurrence of the step  $x^2y^2$ by the two steps $x^2$, $y^2$,
we obtain a path of type B. It is clear that different paths of type A give us different paths af type B 
and each path of type B can be obtained by this operation. The signs of the corresponding monomials 
for this pair of paths are opposite, so each pair contributes zero to the whole ``number'' of paths.

Thus $a_{k,n}$ equals to the number of unpaired paths paths from the point $(0,0)$ to $(k,m)$ with specified steps.
It is clear that the unpaired paths are exactly the paths that do not contain 
neither single steps of the form  $x^2y^2$ nor double steps  $x^2$, $y^2$.
These paths are in an evident one to one correspondence with peakless Motzkin paths.

\paragraph{0-1-2 sums.}

Consider a sum in which the order of summands is fixed, 
each summand is 0, 1 or 2 and for every summand 2 in this sum the next summand is not 0.
We call these sums 0-1-2 \emph{sums}.

Denote by $s_{n,k}$ the number of 0-1-2 sums that are equal to $k$ and consist of $n$ summands
($0\leq k\leq 2n$). For example $s_{3,3}=5$, because the only possible decompositions are the following
$$
3=1+1+1=0+1+2=0+2+1=1+0+2=2+1+0.
$$

\begin{thm}
$s_{n,k}=z_{n,k}$. 
\end{thm}

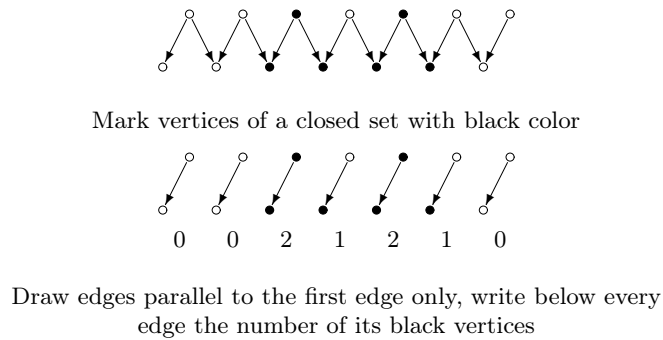
\begin{figure}[h]
\begin{center}\footnotesize
\setlength{\unitlength}{1pt}
%
%
\begin{picture}(130,20)
\multiput(0,0)(20,0){2}{\circle{3}}
\multiput(40,0)(20,0){4}{\circle*{3}}
\put(120,0){\circle{3}}
\multiput(10,20)(20,0){2}{\circle{3}}
\put(50,20){\circle*{3}}\put(90,20){\circle*{3}}
\multiput(110,20)(20,0){2}{\circle{3}}
\put(70,20){\circle{3}}
\multiput(9,18)(20,0){7}{\vector(-1,-2){8}}
\multiput(11,18)(20,0){6}{\vector(1,-2){8}}
\end{picture}

\bigskip
\leavevmode
Mark vertices of a closed set with black color

%
%
\begin{picture}(130,46)(0,-16)
\multiput(0,0)(20,0){2}{\circle{3}}
\multiput(40,0)(20,0){4}{\circle*{3}}
\put(120,0){\circle{3}}
\multiput(10,20)(20,0){2}{\circle{3}}
\put(50,20){\circle*{3}}\put(90,20){\circle*{3}}
\multiput(110,20)(20,0){2}{\circle{3}}
\put(70,20){\circle{3}}
\multiput(9,18)(20,0){7}{\vector(-1,-2){8}}
\put(4,-14){0} \put(24,-14){0} \put(44,-14){2} \put(64,-14){1}
\put(84,-14){2}\put(104,-14){1}\put(124,-14){0}
\end{picture}

\bigskip
\leavevmode
\vbox{\hsize=259pt\noindent 

Draw edges parallel to the first edge only, write below every edge the number of its black vertices}
\end{center}

\kern-3mm
\centering
\caption{
A bijection between closed sets of the graph $Z_{2n}$ and 0-1-2 sums.}
\label{ris:bij2}
\end{figure}

\proof See fig.\ \ref{ris:bij2}. \endproof

\setrisbox\vbox{\hsize=52pt\noindent
\begin{picture}(54,24)
\multiput(0,-2)(10,0){6}{\line(0,1){24}}
\multiput(-2,0)(0,10){3}{\line(1,0){54}}
\Thicklines
\put(0,10){\vector(1,0){10}}
\put(20,10){\vector(1,1){10}}
\put(40,10){\vector(1,-1){10}}
\end{picture}
}

\medskip
We present one more bijection for 0-1-2 sums.

\begin{thm}
$m_{n,k-n}=s_{n,k}$.
\end{thm}

\def\UU{\smash[b]{\lower3.5pt\vbox{\hsize=16pt\noindent 
\begin{picture}(14,14)(-2,-2)
\multiput(0,-2)(10,0){2}{\line(0,1){14}}
\multiput(-2,0)(0,10){2}{\line(1,0){14}}
\Thicklines
\put(0,0){\vector(1,1){10}}
\end{picture}}}}

\def\HH{\lower3.5pt\vbox{\hsize=16pt\noindent 
\begin{picture}(14,14)(-2,-2)
\multiput(0,-2)(10,0){2}{\line(0,1){14}}
\multiput(-2,0)(0,10){2}{\line(1,0){14}}
\Thicklines
\put(0,0){\vector(1,0){10}}
\end{picture}}}

\def\DD{\smash[b]{\lower3.5pt\vbox{\hsize=16pt\noindent 
\begin{picture}(14,14)(-2,-2)
\multiput(0,-2)(10,0){2}{\line(0,1){14}}
\multiput(-2,0)(0,10){2}{\line(1,0){14}}
\Thicklines
\put(0,10){\vector(1,-1){10}}
\end{picture}}}}

\def\UUDD{\lower3.5pt\vbox{\hsize=23pt\noindent 
\begin{picture}(24,14)(-2,-2)
\multiput(0,-2)(10,0){3}{\line(0,1){14}}
\multiput(-2,0)(0,10){2}{\line(1,0){24}}
\Thicklines
\put(0,0){\vector(1,1){10}}
\put(10,10){\vector(1,-1){10}}
\end{picture}}}

\proof Observe that the translation along the vector \UU changes the sum of coordinates by 2,
the translation along the vector \HH changes the sum of coordinates by 1
and the translation along \DD does not change the sum of coordinates.
Consider a peakless Motzkin path from $(0,0)$ to $(n,k-n)$.
It consists of $n$ steps and  translation along this path change the sum of coordinates by  $k$.
Let us construct a sum: in this path replace each segment \UU by 2, each segment \HH by 1, and each segment \DD by 0.
After that put plus signs between the numbers.
Since the path does not contain peaks \UUDD{} , no 0 follows 2. Therefore we obtain a 0-1-2 sum that is equal to $k$ and
consists of $n$ summands. The inverse map is evident so this is a bijection.
Thus $m_{n,k-n}=s_{n,k}$.
\endproof

\paragraph{Weighted paths.}

Let us depict Motzkin paths and claim a payment for drawing of each segment.
A~sloped rising and falling segments cost \$\,1.5,
and we suggest two kind of horizontal segments: ``cheap'' for \$\,1 and
``luxury'' for \$\,2.

\begin{center}
\begin{picture}(110,14)(-5,-11)
\linethickness{0.25pt}
\multiput(0,-2)(10,0){11}{\line(0,1){14}}
\multiput(-2,0)(0,10){2}{\line(1,0){104}}
\Thicklines
\put(0,0){\line(1,0){10}}
\put(60,0){\line(1,1){10}}
\put(90,10){\line(1,-1){10}}
\dottedline[\line(1,1){2}]{2}(30,0)(40,0)
{\footnotesize
\put(-2,-10){\$\,1}
\put(29,-10){\$\,2}
\put(56,-10){\$\,1,5}
\put(87,-10){\$\,1,5}}
\end{picture}
\end{center}

\noindent
Denote by $r_k$ the number of paths that cost \$\,$k$ and such that
their first and last points belong to the same horizontal.
For example  $r_3=5$ as we can see in the following picture.

\begin{center}
\begin{picture}(248,22)(-3,-2)
\linethickness{0.25pt}
\matrixput(0,-2)(10,0){4}(55,0){5}{\line(0,1){24}}
\matrixput(-2,0)(0,10){3}(55,0){5}{\line(1,0){34}}
\Thicklines
\path(0,10)(30,10)
\path(55,10)(65,20)(75,10)
\path(110,10)(120,0)(130,10)
\path(165,10)(175,10)
\dottedline[\line(1,1){2}]{2}(175,10)(185,10)
\path(230,10)(240,10)
\dottedline[\line(1,1){2}]{2}(220,10)(230,10)
\end{picture}
\end{center}

\goodbreak

\begin{thm}
$r_k=a_{k,k}$. 
\end{thm}

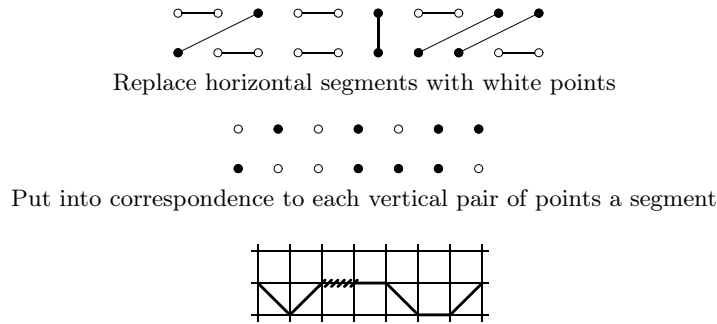
\begin{figure}[h]
\begin{center}\footnotesize
\setlength{\unitlength}{1.5pt}
%
%
\begin{picture}(95,10)
\multiput(0,0)(10,0){10}{\circle{2}}
\multiput(0,10)(10,0){10}{\circle{2}}
\put(0,0){\circle*{2}}
\put(11,0){\line(1,0){8}}
\put(31,0){\line(1,0){8}}
\put(50,0){\circle*{2}}
\put(60,0){\circle*{2}}
\put(70,0){\circle*{2}}
\put(81,0){\line(1,0){8}}
\put(1,10){\line(1,0){8}}
\put(20,10){\circle*{2}}
\put(31,10){\line(1,0){8}}
\put(50,10){\circle*{2}}
\put(61,10){\line(1,0){8}}
\put(80,10){\circle*{2}}
\put(90,10){\circle*{2}}
\put(0,0){\line(2,1){20}}
\put(50,0){\line(0,1){10}}
\put(60,0){\line(2,1){20}}
\put(70,0){\line(2,1){20}}
\end{picture}

\smallskip
\leavevmode
\vbox{\hsize=229pt\noindent 
Replace horizontal segments with white points}
\bigskip

%
%
\begin{picture}(65,10)
\multiput(0,0)(10,0){7}{\circle{2}}
\multiput(0,10)(10,0){7}{\circle{2}}
\put(0,0){\circle*{2}}
\put(30,0){\circle*{2}}
\put(40,0){\circle*{2}}
\put(50,0){\circle*{2}}
\put(10,10){\circle*{2}}
\put(30,10){\circle*{2}}
\put(50,10){\circle*{2}}
\put(60,10){\circle*{2}}
\end{picture}

\smallskip
\leavevmode
\vbox{\hsize=279pt\noindent
Put into correspondence to each vertical pair of points a segment}
\bigskip

%
%
\setlength{\unitlength}{1.2pt}
\begin{picture}(75,24)(-3,-2)
\linethickness{0.25pt}
\multiput(0,-2)(10,0){8}{\line(0,1){24}}
\multiput(-2,0)(0,10){3}{\line(1,0){74}}
\Thicklines
\path(0,10)(10,0)(20,10)
\dottedline[\line(1,1){2}]{2}(20,10)(30,10)
\path(30,10)(40,10)(50,0)(60,0)(70,10)
\end{picture}
\end{center}

\kern-3mm

\centering
\caption{A bijection between partitions of point into pairs and weighted paths}
\label{pic:bij3}
\end{figure}

\proof
In each partition of points into pairs remove  segments that connect points 
from different lines. We obtain a configuration of points and horizontal segments on two lines.
Every ``free'' point color in black, every horizontal segment replace by one white point and 
place all these points on the lines uniformly.
Then put into correspondence to each vertical pair of points a segment of a path (see fig.~\ref{pic:bij3}):
let a pair \
\begin{picture}(2,7)
\put(0,0){\circle*{2.5}}
\put(0,7){\circle*{2.5}}
\end{picture}
correspond to cheap horizontal segment, a pair \
\begin{picture}(2,7)
\put(0,0){\circle{2.5}}
\put(0,7){\circle{2.5}}
\end{picture}
correspond to luxury horizontal segment, a pair \
\begin{picture}(2,7)
\put(0,0){\circle{2.5}}
\put(0,7){\circle*{2.5}}
\end{picture}
correspond to segment \UU, and a pair \
\begin{picture}(2,7)
\put(0,0){\circle*{2.5}}
\put(0,7){\circle{2.5}}
\end{picture}
correspond to segment 
\DD. \endproof

\paragraph{Symmetric configurations of chords.}

Fix an arbitrary integers $\ell>2$ and $n$.

Consider an inscribed $\ell n$-gon. Split its circle onto $\ell$ equal arcs
whose endpoints do not coincide with vertices of the polygon. 
Join some vertices by chords.
We are interested in configurations of chords that satisfy the following properties:

1) a configuration has symmetry of $\ell$-th order, i.e. it is invariant under rotations by angle  $2\pi/\ell$;

2) chords do not intersect (and have no common endpoints); 

3) endpoints of any chord can not be neighboring points on the same arc.

\begin{figure}[h]
\begin{center}
\setlength{\unitlength}{1.2pt}\footnotesize
\leavevmode
\qquad\qquad\qquad\hfill\epsfbox{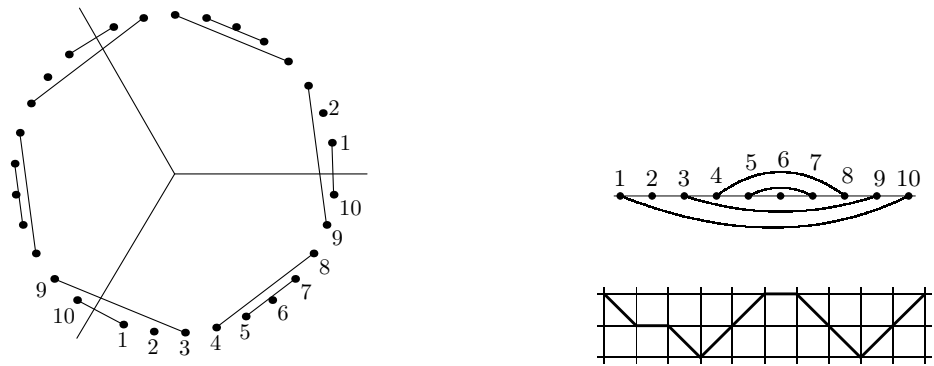}
\hfill
\vbox{\hsize=10cm
\begin{picture}(94,35)(-2,-20)
\path(-2,0)(92,0)
\multiput(0,0)(10,0){10}{\circle*{2}}
\qbezier(0,0)(50,-20)(90,0)
\qbezier(20,0)(50,-10)(80,0)
\qbezier(30,0)(50,15)(70,0)
\qbezier(40,0)(50,5)(60,0)
\put(-2,3){1}   \put(8,3){2}   \put(18,3){3}
\put(28,4){4}  \put(39,8){5}   \put(49,9){6}
\put(59,8){7}  \put(69,4){8}   \put(79,3){9}
\put(86,3){10}
\end{picture}

\bigskip
\begin{picture}(100,20)
\linethickness{0.25pt}
\multiput(0,-2)(10,0){11}{\line(0,1){24}}
\multiput(-2,0)(0,10){3}{\line(1,0){104}}
\Thicklines
\path(0,20)(10,10)(20,10)(30,0)(40,10)(50,20)(60,20)(70,10)(80,0)(90,10)(100,20)
\end{picture}}
\end{center}

\kern-2mm

\centering
\caption{
A configuration with symmetry of 3rd order, its arc diagram and peakless Motzkin path.
}
\label{pic:chords}

\end{figure}

It follows that if endpoints of a chord are on different arcs then these arcs are neighboring.
So the configuration is determined by a picture that we see in any of $\ell$ sectors,
and the number of configurations does not depend on $\ell$.
An example of a configuration for $n=10$, $\ell=3$ is shown in fig.~\ref{pic:chords}.

\begin{thm}
The number of symmetric chord configurations equals $a_{n,n}$.
\end{thm}

\proof 
Choose one of the $\ell$ arcs and ``straighten'' it.
We obtain a segment with $n$ marked points. 
Depict a chord that joins two points in the initial sector as an arc above the segment.
Depict a chord that join a point from the initial sector, say the point number $p$, 
and a point from the neighboring sector, say the point number $q$ in that sector, as an arc 
below the segment (see fig.~\ref{pic:chords}).
Now we can easily transform this arc diagram into peakless Motzkin path of length $n$.
To do this we inspect all the marked points from left to right.
If we see an arc that rises from the point or we see an an arc goes to the point from below we depict a rising segment of path;
isolated points we depict a horizontal edges, 
and if an arc goes to the point from above or starts from the point and goes below the segment we depict a falling segment of 
Motzkin path.
\endproof

\epsfxsize=1.9cm
\setrisbox\vbox{\hsize=2cm\noindent\epsfbox{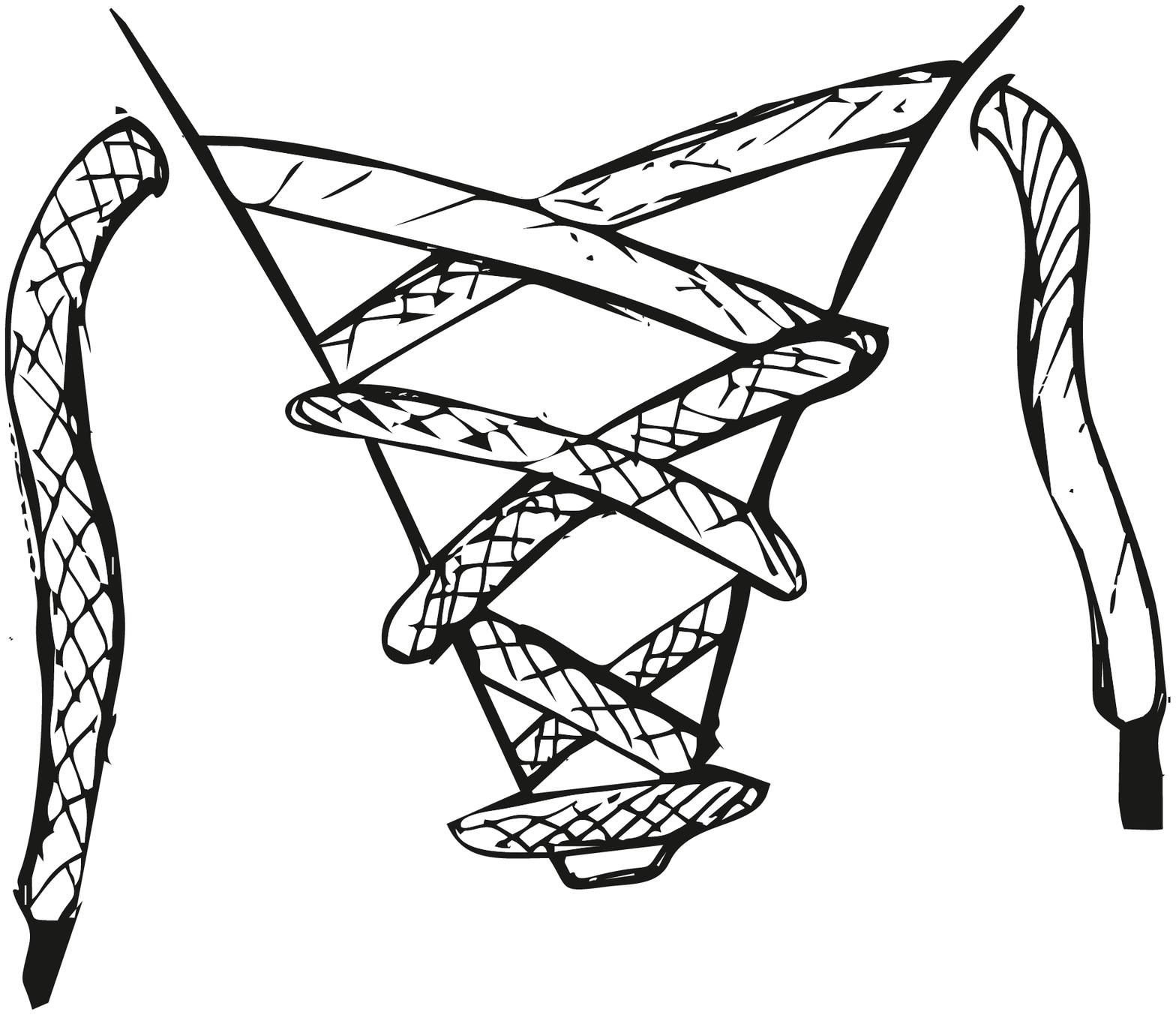}}

\paragraph{Lacing.} 

\ris 
Let us count  the number of ways to lace a shoe that has two lines with $n$ holes in each for lacing.
We assume that the following restrictions hold:
\endris

1) the lacing starts and finishes at the topmost pair of holes;

2) the shoelace passes through each hole exactly one time;

3) for every hole at least one of its neighboring holes along the shoelace is at the opposite line of holes
(the first and the last hole will satisfy this property after a knot has been tied);

4) we do not take into account topological details like knotting of shoelace, is a shoelace passes a hole upwards or top-down, etc.

\begin{figure}[h]
\kern-2mm
\begin{center}
\setlength{\unitlength}{1.2pt}
\begin{picture}(20,50)
\multiput(0,0)(0,15){4}{\circle{3}}
\multiput(15,0)(0,15){4}{\circle{3}}
\put(0,45){\line(-1,1){6}}
\put(15,45){\line(1,1){6}}
\put(0,30){\line(0,1){15}}
\put(0,30){\line(1,-2){15}}
\put(0,0){\line(1,0){15}}
\put(0,0){\line(1,1){15}}
\put(15,15){\line(0,1){15}}
\put(0,15){\line(1,1){15}}
\put(0,15){\line(1,2){15}}
\end{picture}
\qquad
\begin{picture}(20,50)
\multiput(0,0)(0,15){4}{\circle{3}}
\multiput(15,0)(0,15){4}{\circle{3}}
\put(0,45){\line(-1,1){6}}
\put(15,45){\line(1,1){6}}
\put(0,30){\line(0,1){15}}
\put(15,30){\line(0,1){15}}
\put(15,15){\line(-1,1){15}}
\put(0,15){\line(1,1){15}}
\put(15,0){\line(0,1){15}}
\put(0,0){\line(0,1){15}}
\put(0,0){\line(1,0){15}}
\end{picture}
\qquad\qquad\qquad\qquad
\begin{picture}(20,50)
\multiput(0,0)(0,15){4}{\circle{3}}
\multiput(15,0)(0,15){4}{\circle{3}}
\put(0,45){\line(-1,1){6}}
\put(15,45){\line(1,1){6}}
\put(0,30){\line(0,1){15}}
\put(0,15){\line(0,1){15}}
\put(15,30){\line(0,1){15}}
\put(15,15){\line(0,1){15}}
\put(0,15){\line(1,-1){15}}
\put(0,0){\line(1,1){15}}
\put(0,0){\line(1,0){15}}
\end{picture}
\qquad
\begin{picture}(20,50)
\multiput(0,0)(0,15){4}{\circle{3}}
\multiput(15,0)(0,15){4}{\circle{3}}
\put(0,45){\line(-1,1){6}}
\put(15,30){\line(1,1){6}}
\put(0,30){\line(0,1){15}}
\put(0,30){\line(1,-2){15}}
\put(0,0){\line(1,0){15}}
\put(0,0){\line(1,1){15}}
\put(15,30){\line(0,1){15}}
\put(0,15){\line(1,0){15}}
\put(0,15){\line(1,2){15}}
\end{picture}
\end{center}

\kern-4mm
\caption{Two right (on the left) and two ``wrong'' (on the right) lacings.} 
\label{ris:lacing}
\end{figure}
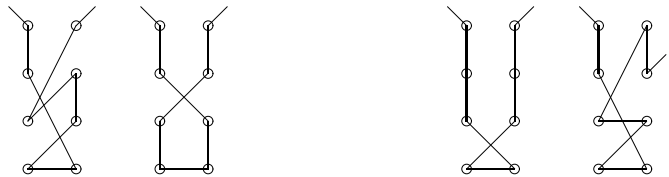

We call a lacing that satisfies these restrictions \emph{right} (see fig.~\ref{ris:lacing}).
We call a lacing \emph{non self-crossing} if it starts at the upper left hole and finishes at the lower right hole,
has no self-crossings and satisfies properties 2)--4).

\begin{thm}
The number of non self-crossing lacings  equals $a_{n,n}$.
The number of right lacings  equals $((n-1)!)^2a_{n,n}$.
\end{thm}

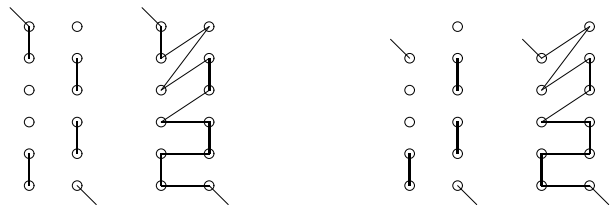
\begin{figure}[b]
\kern-2mm
\begin{center}
\setlength{\unitlength}{1.2pt}
\begin{picture}(20,55)
\multiput(0,0)(0,10){6}{\circle{3}}
\multiput(15,0)(0,10){6}{\circle{3}}
\put(0,50){\line(-1,1){6}}
\put(15,0){\line(1,-1){6}}
\put(0,40){\line(0,1){10}}
\put(0,0){\line(0,1){10}}
\put(15,10){\line(0,1){10}}
\put(15,30){\line(0,1){10}}
\end{picture}
\qquad
\begin{picture}(20,55)
\multiput(0,0)(0,10){6}{\circle{3}}
\multiput(15,0)(0,10){6}{\circle{3}}
\put(0,50){\line(-1,1){6}}
\put(15,0){\line(1,-1){6}}
\put(0,40){\line(0,1){10}}
\put(0,0){\line(0,1){10}}
\put(15,10){\line(0,1){10}}
\put(15,30){\line(0,1){10}}
\put(0,40){\line(3,2){15}}
\put(0,30){\line(3,4){15}}
\put(0,30){\line(3,2){15}}
\put(0,20){\line(3,2){15}}
\put(0,20){\line(1,0){15}}
\put(0,10){\line(1,0){15}}
\put(0,0){\line(1,0){15}}
\end{picture}
\qquad\qquad\qquad
%
%
\begin{picture}(20,55)
\multiput(0,0)(0,10){5}{\circle{3}}
\multiput(15,0)(0,10){6}{\circle{3}}
\put(0,40){\line(-1,1){6}}
\put(15,0){\line(1,-1){6}}
\put(0,0){\line(0,1){10}}
\put(15,10){\line(0,1){10}}
\put(15,30){\line(0,1){10}}
\end{picture}
\qquad
\begin{picture}(20,55)
\multiput(0,0)(0,10){5}{\circle{3}}
\multiput(15,0)(0,10){6}{\circle{3}}
\put(0,40){\line(-1,1){6}}
\put(15,0){\line(1,-1){6}}
\put(0,0){\line(0,1){10}}
\put(15,10){\line(0,1){10}}
\put(15,30){\line(0,1){10}}
\put(0,40){\line(3,2){15}}
\put(0,30){\line(3,4){15}}
\put(0,30){\line(3,2){15}}
\put(0,20){\line(3,2){15}}
\put(0,20){\line(1,0){15}}
\put(0,10){\line(1,0){15}}
\put(0,0){\line(1,0){15}}
\end{picture}
\end{center}

\kern-2mm
\caption{\vtop{\hsize=307pt 
A configuration of points and segments determines uniquely a non self-crossing lacing,
even on a defective shoe.}}
\label{ris:vertik-uchastki}

\kern-5mm
\end{figure}

\proof
First observe that a non self-crossing lacing is uniquely determined by the set of its vertical segments 
connecting neighboring vertical holes (on the same side of the shoe).
Indeed, if a set of vertical segments then the topmost ``free'' holes on the opposite sides of the shoe should be neighboring in our lacing, 
(see fig.~\ref{ris:vertik-uchastki}, left), after that the next two ``free'' holes on the opposite sides should be neighboring too, etc.
Therefore the number of non self-crossing lacings equals the number of ways to choose a set vertical segments, i.e. $a_{n,n}$.

Given a non self-crossing lacing we can construct $(n!)^2$ new lacings by 
permuting holes in each side. The lacings obtained in this way can start and finish 
in arbitrary holes on the opposite sides of the shoe.
Conversely, an arbitrary lacing up to a permutation of holes in each side determines uniquely
a non self-crossing lacing: it is a lacing obtained by the permutation for which the movement along the lace
always deliver us to the topmost free holes (on the corresponding side of the shoe).

Thus the number of lacings without restriction given by the property 1) equals  $(n!)^2a_{n,n}$.
Right lacings necessarily start and finish at the two topmost holes, therefore the number 
of right lacings equals $((n-1)!)^2a_{n,n}$.
\endproof

\paragraph{Defective shoe.}

We call a shoe \emph{defective} if it has different number of holes on its sides: 
$k$ holes on the left side and $n$ holes on the right side.
Consider lacings of defective shoe. As in previous paragraph we call a lacing \emph{right} 
if it satisfies properties 1)--4). The definition of 
\emph{non self-crossing} lacing is word by word the same as in the previous paragraph, too.

\goodbreak


Again, each non self-crossing lacing determines a configuration of 
points and segments on the sides of the shoe
(see fig.~\ref{ris:vertik-uchastki}, right).
Each element of the configuration --- a point or a segment --- we will call an \emph{object}.
Each object has two neighboring objects (along the lace), both of them are on the opposite side of the shoe.
Therefore each configuration that is obtained from some lacing has equal number of objects on both sides.
Conversely, each configuration with equal number of objects on both sides determines uniquely some non self-crossing lacing.


\medskip

\begin{figure}[h]
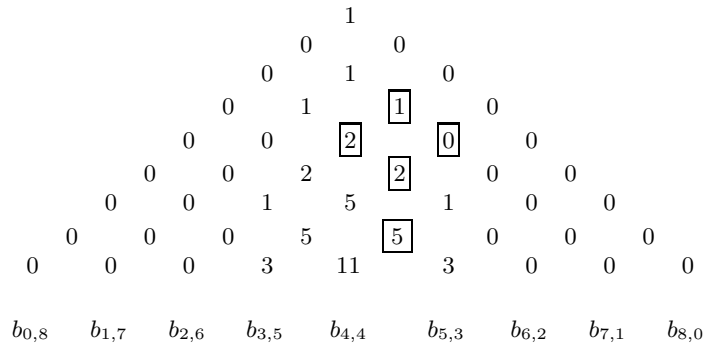

{\footnotesize
$$\setcounter{MaxMatrixCols}{20}
\begin{matrix}\setlength{\arraycolsep}{0pt}
 & & & & & & & &1\\
 & & & & & & &0& &0\\
 & & & & & &0& &1& &0\\
 & & & & &0& &1& &\boxed{\!1\!}&&0\\
 & & & &0& &0& &\boxed{\!2\!}&&\boxed{\!0\!}&&0\\
 & & &0& &0& &2& &\boxed{\!2\!}&&0&&0\\
 & &0& &0& &1& &5& &1& &0& &0\\
 &0& &0& &0& &5& &\!\!\!\!\boxed{5}\!\!\!\!& &0& &0& &0\\ 
0& &0& &0& &3& &\!\!11\!\!&&3& &0& &0& &0\\[13pt]
\!\!\!\!b_{0,8}\!\!\!\! &&\!\!\!\!b_{1,7}\!\!\!\! && \!\!\!\!b_{2,6}\!\!\!\!  &&  \!\!\!\!b_{3,5}\!\!\!\! &&  \!\!\!\!b_{4,4}\!\!\!\!&&
\!\!\!\!b_{5,3}\!\!\!\! && \!\!\!\!b_{6,2}\!\!\!\!  &&  \!\!\!\!b_{7,1}\!\!\!\! &&  \!\!\!\!b_{8,0}\!\!\!\!              \\
\end{matrix}$$
\caption{A triangle $b_{k,n}$ and example of its recurrence relation.}
\label{ris:b_triangle}

\kern-5pt}
\end{figure}

In this connection consider one more combinatorial sequence.
Let two parallel lines are given, $k$ points are marked on the first line, $n$ points are marked on the second line and
some pairs of neighboring points on the lines are joined by segments, the segments do not intersect.
Denote by $b_{k,n}$ the number of these configurations which have equal number of objects on the lines.
Define also a ``degenerate'' values of $b_{k,n}$: let  $b_{0,0}=1$ and $b_{k,n}=0$, if $(k,n)\ne (0,0)$, 
but at least one of the numbers $k$, $n$ in non positive.
Values of the sequence $b_{k,n}$ for small $k$ and $n$ are shown in fig.~\ref{ris:b_triangle}.
The topmost 1 is $b_{0,0}$, the numbering of elements in the last row is shown for convenience.

Analogously to $a_{k,n}$, if we remove the leftmost objects in each side of a configuration, 
we obtain a recurrence
$$
b_{k,n}=b_{k-1,n-1}+b_{k-1,n-2}+b_{k-2,n-1}+b_{k-2,n-2}.
$$
This recurrence holds for all $n$ and $k$, except $(k,n)=(0,0)$.
Due to this recurrence we immediately obtain a formula for the generation function:
\begin{equation}
B(x,y)=\sum_{k,n} b_{k,n}  x^k y^n =\frac1{1-(xy+ x^2y +x y^2+x^2y^2)}-1. 
\label{eqn:B-generatin-function}
\end{equation}

\begin{thm}
$a_{n,n}=b_{n,n}$.
\end{thm}

\proof
It follows from the definition of $b_{k,n}$.
\endproof

By reasoning like in the previous paragraph we conclude that 
the sequence  $b_{k,n}$ allows us to count lacings:
the number of non-crossing lacings is $b_{k,n}$,
the number of arbitrary lacings is $ (k-1)!(n-1)! b_{k,n}$,
and the number of right lasings is $ k!\,n!\,b_{k,n}$.
As a corollary we obtain that the function~\eqref{eqn:B-generatin-function}
is an exponential generation function of the numbers of arbitrary lacings.


\medskip

Let us describe several more combinatorial realizations of the sequence $b_{k,n}$.

Since the generating function \eqref{eqn:B-generatin-function} of the sequence $b_{k,n}$
is similar to the generation function of the sequence $a_{k,n}$, we can 
apply reasoning of the example 2 and obtain that $b_{k,n}$ is a number of paths 
on the plane from the origin to the point $(k,n)$ with  segments of four types depicted below.

\begin{center}
\begin{picture}(24,34)(-2,-12)
\multiput(0,-2)(10,0){3}{\line(0,1){24}}
\multiput(-2,0)(0,10){3}{\line(1,0){24}}
\Thicklines
\put(0,0){\vector(1,1){10}}
\put(5,-12){$xy$}
\end{picture}
\qquad
\begin{picture}(24,34)(-2,-12)
\multiput(0,-2)(10,0){3}{\line(0,1){24}}
\multiput(-2,0)(0,10){3}{\line(1,0){24}}
\Thicklines
\put(0,0){\vector(1,1){20}}
\put(2,-12){$x^2y^2$}
\end{picture}
\qquad
\begin{picture}(24,34)(-2,-12)
\multiput(0,-2)(10,0){3}{\line(0,1){24}}
\multiput(-2,0)(0,10){3}{\line(1,0){24}}
\Thicklines
\put(0,0){\vector(2,1){20}}
\put(4,-12){$x^2y$}
\end{picture}
\qquad
\begin{picture}(24,34)(-2,-12)
\multiput(0,-2)(10,0){3}{\line(0,1){24}}
\multiput(-2,0)(0,10){3}{\line(1,0){24}}
\Thicklines
\put(0,0){\vector(1,2){10}}
\put(4,-12){$xy^2$}
\end{picture}
\end{center}

\goodbreak

Now, rotating a picture by  $45^\circ$ (and taking a diagonal of a square as unity)
we see that $b_{k,n}$ equals the number of weighted Motzkin paths from the origin to
the point $(\frac{k+n}{2},\frac{k-n}{2})$. The prices are exactly the lengths of
projections onto the line  $y=x$.
Thus we have check once again that $b_{k,k}=r_k$.



\medskip
\begin{center}
\begin{picture}(110,14)(-5,-11)
\linethickness{0.25pt}
\multiput(0,-2)(10,0){11}{\line(0,1){14}}
\multiput(-2,0)(0,10){2}{\line(1,0){104}}
\Thicklines
\put(0,0){\line(1,0){10}}
\put(60,0){\line(1,1){10}}
\put(90,10){\line(1,-1){10}}
\dottedline[\line(1,1){2}]{2}(30,0)(40,0)
{\footnotesize
\put(-2,-10){\$\,1}
\put(29,-10){\$\,2}
\put(56,-10){\$\,1,5}
\put(87,-10){\$\,1,5}}
\end{picture}
\end{center}


Finally, observe that by definition $b_{k,n}$ is equal to the number of pairs ``$S_1$-composition 
of $k$ and $S_1$-composition of $n$'' with equal numbers of summands.
We can depict these pairs of compositions as ``staircases'' from (0,0) to ($k$, $n$).
A \emph{staircase} is a path in which vertical and horizontal segments alternate,
the first segment is horizontal, the last segment is vertical, and the lengths of all segments are 1 or 2. 
For example a pair of compositions $8=2+2+1+1+2=2+1+2+1+2$ that we have seen in fig.~\ref{ris:bij1},
determines the following staircase.

\begin{center}
\begin{picture}(80,80)
\linethickness{0.25pt}
\multiput(0,0)(10,0){9}{\line(0,1){80}}
\multiput(0,0)(0,10){9}{\line(1,0){80}}
\Thicklines
\path(0,0)(20,0)(20,20)(40,20)(40,30)(50,30)(50,50)(60,50)(60,60)(80,60)(80,80)
{\footnotesize
\put(-19,-7){(0,0)}\put(83,78){(8,8)}
\put(8,-10){2}\put(28,-10){2}\put(43,-10){1}\put(53,-10){1}\put(68,-10){2}
\put(-8,8){2}\put(-8,22){1}\put(-8,38){2}\put(-8,52){1}\put(-8,68){2}}
\end{picture}
\end{center}

\paragraph{Diagonal.} 
In this paragraph we use theory of functions of a complex variable.
The reader, who has not studied this course yet, can see details in \cite{Shabat}.

We call a sequence $a_{n,n}$ a \emph{diagonal} sequence.
We start form the formula for the generation function of the sequence  $a_{n,k}$:
$$
F(x,y)=\sum_{k,n} a_{k,n}  x^k y^n=\frac1{1-(x^2+y^2+xy-x^2y^2)}.
$$
Due to the estimation \eqref{eqn:ocenka_a_nk_sverhu}, this series converges for $|x|<\varphi^{-1}$, $|y|<\varphi^{-1}$,
where $\varphi=\frac{1+\sqrt{5}}2$ is a golden ratio.

Now find the generating function of diagonal sequence:
$g(x)=\smash[t]{\sum\limits_{k=0}^{+\infty} a_{k,k} x^k}$.
We apply the following standard trick.
Fix a sufficiently small $x$ (we may consider $x$ as a real variable) and consider a function
$$
H(s)=F\Bigl(\sqrt{s}, \frac{x}{\sqrt{s}}\Bigr)=\frac{-s}{s^2-s(x^2-x+1)+x^2}.
$$
It can be decomposed in a Laurent series
$\sum\limits_{k,n} a_{k,n}  (\sqrt{s})^k \bigl(\frac x{\sqrt{s}}\bigr)^n$ 
over powers of~$s$ (and non negative powers of $x$).
This series obviously converges on the annulus $|x|^2\varphi^2 < |s|<\varphi^{-2}$, 
and the function $H(s)$ is rational and defined at the whole complex plane.
It is easy to see that $g(x)$ is a constant term in this series.
We can find it by the residue theorem:
$$
g(x)=\frac1{2\pi i}\int_{|s|=\rho} \frac{H(s)\,ds}{s} = \sum\text{Res} \frac{H(s)}{s},
$$
where the integration is over an arbitrary circle inside the annulus
and the the sum runs over the singularities inside the circle.

The discriminant of the denominator of $H(s)$ equals
$$
(x^2-x+1)^2-4x^2=1-2x-x^2-2x^3+x^4.
$$
This expression is approximately equal to 1 for small $x$, 
therefore one of the singularities of $H(s)$ is not far from 1 (outside the integration contour),
an the second singularity 
$$
s_0= \frac{1-x+x^2-\sqrt{1-2x-x^2-2x^3+x^4}}{2}.
$$
is near 0 inside the contour. The residue in this point equals (standard exercise)
$$
g(x)=
\frac1{\sqrt{1-2x-x^2-2x^3+x^4}}.
$$

\paragraph{Recurrence for the diagonal sequence 
$r_n=a_{n,n}$.} 
For generating function $g(x)$ we have
$$
g'(x) = \smash[b]{\frac{1+x+3x^2-2x^3}{(1-2x-x^2-2x^3+x^4)^{3/2}} = g(x)\cdot\frac{1+x+3x^2-2x^3}{1-2x-x^2-2x^3+x^4}.}
$$
Hence 
$$
(1-2x-x^2-2x^3+x^4)g'(x) = (1+x+3x^2-2x^3)g(x) .
$$
Since $g'(x)=\sum\limits_{k=0}^{+\infty} kr_k x^{k-1}$,
equating  coefficients of  $x^{n-1}$ we obtain 
$$
nr_n-2(n-1)r_{n-1}-(n-2)r_{n-2}-2(n-3)r_{n-3}+(n-4)r_{n-4} =
r_{n-1}+r_{n-2}+3r_{n-3}-2r_{n-4} .
$$
Thus, 
$$
nr_n-(2n-1)r_{n-1}-(n-1)r_{n-2}-(2n-3)r_{n-3}+(n-2)r_{n-4} = 0 .
$$

\paragraph{Asymptotics of the diagonals sequence.} 

\def\phi{\varphi}
Applying Darboux method (see~\cite[\S4.3]{Green-Knuth}) we fill find an asymptotic formula for  $r_n$.

The idea of Darboux method is the following. 
If a function is defined as a sum of a power series $g(z)=\sum r_n z^n$, then this series is 
automatically its Taylor expansion at the point $z=0$. Its radius of convergence is equal to the distance from 0
to the nearest singularity point $z_1$ of the function $g(z)$.
If we choose a ``simple'' function $h(z)=\sum h_n z^n$, that has ``the same'' singularity in the point $z_1$,
then it is possible that the difference $g(z)-h(z)=\sum (r_n-h_n) z^n$ is regular in the point $z_1$.
In this case the radius of convergence of $g(z)-h(z)=\sum (r_n-h_n) z^n$ is greater then the radius of $g(z)$, 
that means that the coefficients of the second series is much smaller than the coefficients of the first one,
i.e. $r_n-h_n=o(r_n)$ for $n\to+\infty$. Thus $r_n\sim h_n$.

For the function $g(z)$ the nearest to 0 singular point is a branch point, it can not be canceled just by a subtraction of some function.
So we need a bit more accurate reasoning.

Factor the expression under the square root in the formula for $g(z)$:
$$
1-2z-z^2-2z^3+z^4=(1-\phi^2z)(1-\phi^{-2}z)(1+z+z^2),
$$
where $\phi = \frac{1+\sqrt{5}}{2}$ is a golden ratio. Let $f(z) = \frac{1}{\sqrt{1-\phi^{-2}z}\sqrt{1+z+z^2}}$.
Then 
$g(z)$ is equivalent $\frac{f(\phi^{-2})}{\sqrt{1-\phi^2z}}$ for $z\to \phi^{-2}$
and the difference of these functions is bounded on the neighborhood of the point $z=\phi^{-2}$:
\begin{equation}
g(z)-\frac{f(\phi^{-2})}{\sqrt{1-\phi^2z}}= \frac{f(z)-f(\phi^{-2})}{\sqrt{1-\phi^2z}}=\sqrt{1-\phi^2z} h(z) ,
\label{eqn:g(z)-Darboux}
\end{equation}
where the function $h(z)$ has singular points $\phi^2$ and $e^{\pm 2\pi i/3}$. 

\lemma.
Let $\displaystyle
(1-t)^p = \smash[b]{\sum_{n=0}^\infty} \lambda_{p,n} t^n 
$
be a Taylor series of the function $(1-t)^p$ at the point $t=0$. Then
$$
\lambda_{-1/2,n} = \frac{1}{\sqrt{\pi n}}+O\Big(\frac{1}{n^{3/2}}\Big)
\quad\text{and}\quad
\lambda_{1/2,n} = O\Big(\frac{1}{n^{3/2}}\Big) .
$$

For the proof use Newton's generalized binomial theorem and Stirling's formula.

\begin{thm}
$r_n = \dfrac{\phi^{2n+2}}{2\sqrt[4]{5}\sqrt{\pi n}} + 
O\Big(\dfrac{\phi^{2n}}{n^{3/2}}\Big)$, where $\phi = \frac{1+\sqrt{5}}{2}$.
\end{thm}

\proof 
By  \eqref{eqn:g(z)-Darboux} the number $r_n$ equals the sum of coefficients of  $z^n$ 
in the expansions of functions $\frac{f(\phi^{-2})}{\sqrt{1-\phi^2z}}$ and $\sqrt{1-\phi^2z} h(z)$.

\goodbreak

Consider the coefficient of $z^n$ in the expansion of $\sqrt{1-\phi^2z} h(z)$.
Since the singular points of  $h(z)$ are $\phi^2$ and $e^{\pm 2\pi i/3}$
its Taylor series $h(z)=\sum\alpha_n z^n$ at the point $z=0$ converges in 
any circle of raduis $\rho<1$ with the center in the point $z=0$.
It follows that $|\alpha_n|< c/\rho^n$ 
(otherwise for $z=\frac{\rho+1}{2}$ terms of the series tends to infinity and the series diverges).
So the the coefficient of $z^n$ in the expansion of $\sqrt{1-\phi^2z} h(z)$ equals
\begin{align*}
\sum_{k=0}^n \phi^{2k}\lambda_{1/2,k}\alpha_{n-k} &=
 \alpha_n + \smash[t]{\sum_{k=1}^{[n/2]}} \ldots + \sum_{k=[n/2]+1}^n \ldots = 
\\
&= O\Big(\frac{1}{\rho^n}\Big) + 
  O\Big(\frac{1}{\rho^{n/2}}\Big)\sum_{k=1}^{[n/2]} \phi^{2k}\lambda_{1/2,k} + 
  O\Big(\frac{\phi^{2n}}{n^{3/2}}\Big) \!\!\!\sum_{k=[n/2]+1}^n\!\!\!\alpha_{n-k} = 
\\
&= O\Big(\frac{1}{\rho^n}\Big) + 
  O\Big(\frac{\phi^n}{\rho^{n/2}}\Big) O\Big(\frac{1}{\sqrt{n}}\Big) +
  O\Big(\frac{\phi^{2n}}{n^{3/2}}\Big) O\Big(\frac{1}{1-\rho}\Big) = 
  O\Big(\frac{\phi^{2n}}{n^{3/2}}\Big) .
\end{align*}
Here the signs ``$O$'' with semi-integer powers are obtained by applying the lemma, 
and in the last equality we estimate from above the sum of $\alpha_{n-k}$ by geometric series.

To finish the proof it remains to observe that the coefficient of $z^n$ in the expansion of
$\smash{\frac{f(\phi^{-2})}{\sqrt{1-\phi^2z}}}$ is much greater, it equals
$$
f(\phi^{-2})\phi^{2n}\lambda_{-1/2,n} = 
  \frac{f(\phi^{-2})\phi^{2n}}{\sqrt{\pi n}}+O\Big(\frac{\phi^{2n}}{n^{3/2}}\Big) 
= \frac{f(\phi^{-2})\phi^{2n}}{\sqrt{\pi n}}+O\Big(\frac{\phi^{2n}}{n^{3/2}}\Big) = 
  \frac{\phi^{2n+2}}{2\sqrt[4]{5}\sqrt{\pi n}}+O\Big(\frac{\phi^{2n}}{n^{3/2}}\Big) .
\qedhere
$$
\endproof

\smallbreak
For any integer $n$ one can obtain in this way an asymptotic expansion of $r_n$ with accuracy 
$O\Big(\dfrac{\phi^{2n}}{n^{m+1/2}}\Big)$ by exploring more terms in 
expansion of $g(z)$ and applying the formula (2.2)~\cite{Flajolet-Odlyzko} instead of our lemma.
The second term of the asymptotic expansion is also found in~\cite[equation~(24)]{Hofacker-Reidys-Stadler}
with different approach.


\paragraph{Sources.}

All the sequences under discussion can be found in the On-Line Encyclopedia of Integer Sequences \cite{OEIS}:
a ``triangular'' sequence  $a_{k,n}$ (or speaking more accurately $z_{n,k}$) is the sequence A079487, 
sequence $b_{k,n}$ is A125250, the diagonal sequence $a_{n,n}$ is A051286 and lacing sequence is  A078698.
The main combinatorial realization of the sequences are also given.
Detailed study of the sequence  $a_{k,n}$ on the language of order ideals of posets is in  \cite{munarini},
we call this realizations ``fences'', the term ``Whitney numbers of the second type'' is also used.
Symmetric configurations of chords are described in~\cite{Hofacker-Reidys-Stadler}.
Weighted paths and compositions can be found in \cite{Bona-Knopfmacher}, 
domino and compositions are in \cite{Banderier}, lacings are in \cite{Ken}.
Simplified version of this article is published in \cite{spb14}.



\end{document}